\documentclass[11pt]{article}
\usepackage[a4paper,margin=2.7cm]{geometry}
\usepackage{amsmath,amssymb,amsthm,mathtools}
\usepackage{mathrsfs}
\usepackage{enumitem}
\usepackage{xcolor}
\usepackage[colorlinks=true,linkcolor=blue,citecolor=blue,urlcolor=blue]{hyperref}
\usepackage{microtype}
\usepackage{bm}
\usepackage{cite}
\usepackage[utf8]{inputenc}
\newtheorem{theorem}{Theorem}[section]
\newtheorem*{theorem*}{Theorem}
\newtheorem{proposition}[theorem]{Proposition}
\newtheorem{lemma}[theorem]{Lemma}
\newtheorem{corollary}[theorem]{Corollary}
\newtheorem{assumption}[theorem]{Assumption}
\newtheorem{definition}[theorem]{Definition}
\theoremstyle{remark}
\newtheorem{remark}[theorem]{Remark}
\newtheorem{example}[theorem]{Example}
\newcommand{\E}{\mathbb E}
\newcommand{\Pp}{\mathbb P}
\newcommand{\R}{\mathbb R}
\newcommand{\Law}{\mathcal L}
\newcommand{\eps}{\varepsilon}
\newcommand{\dd}{\,\mathrm d}
\newcommand{\A}{\mathcal A}
\newcommand{\Lc}{\mathcal L}
\newcommand{\Lav}{\overline{\mathcal L}}
\newcommand{\PN}{\Pi_N}
\newcommand{\nN}{n_N}
\newcommand{\Pol}{\Upsilon}
\newcommand{\HS}{\mathrm{HS}}
\newcommand{\sff}{\mathsf f}
\newcommand{\sfs}{\mathsf b}
\newcommand{\T}{\mathbb T}

\title{From Random Advection to Transport Noise and the LU Framework for the Navier--Stokes Equations}
\author{%
Arnaud Debussche$^{2,3}$ \qquad B\'erenger Hug$^{1}$ \qquad \'Etienne M\'emin$^{1}$\\[0.8em]
\small $^{1}$Univ Rennes, Inria, Odyssey, IRMAR-UMR 6625,\\ \small Centre Inria Rennes 35042 Rennes Cedex, France\\
\small $^{2}$Univ Rennes, CNRS, IRMAR-UMR 6625, F-35000 Rennes, France\\
\small $^{3}$Institut universitaire de France (IUF)%
}
\date{}

\begin{document}
\maketitle

\begin{abstract}
We consider the incompressible Navier--Stokes equations in a smooth bounded domain in dimension two or three, with Dirichlet boundary conditions, when the velocity is advected by a rapidly varying stationary Markov process $\eps^{-1}m(t/\eps^2)$. The fast field may be spatially inhomogeneous. We prove that, when $\eps\to0$, the laws of the large-scale velocities are tight and that their limit points are martingale solutions of the stochastic Navier--Stokes equations with transport noise in the form of the location uncertainty framework. The limiting equation contains in particular a noise which is described in terms of the covariance of the fast advection. Under a symmetry assumption on the temporal correlations of the fast field, the noise is Stratonovich and an advection by the It\^o--Stokes drift appears. The proof relies on perturbed test functions, which we use both for the uniform energy bounds and for the identification of the limit. We have to adapt this method and add the Stokes operator to the Poisson equation associated to the fast advection in order to use fewer regularity assumptions and obtain correctors which satisfy the boundary conditions.
\end{abstract}



\section{Introduction}

In this work, we study the incompressible Navier--Stokes equations advected by a random velocity field which decorrelates rapidly in time, and we consider the diffusion limit when the correlation time goes to zero. This question is motivated by the modeling of unresolved scales in fluid dynamics, where stochastic transport models are now widely used. The idea is to split the Lagrangian velocity into a large-scale component, regular in time, and a rapidly fluctuating one. In the modeling under location uncertainty (LU), the fluctuating part is a velocity field which is white in time and correlated in space. Using this decomposition in the Reynolds transport theorem and in the conservation laws of an incompressible fluid, we obtain a stochastic Navier--Stokes equation with a transport noise, an effective diffusion, and a drift due to the spatial inhomogeneity of the unresolved velocity. We refer to M\'emin~\cite{Memin2014} and Resseguier, M\'emin and Chapron~\cite{ResseguierMeminChapron2017} for the derivation of LU equations in fluid dynamics, and to~\cite{ResseguierMeminHeitzChapron2017,BauerChandramouliChapronLiMemin2020} for applications, in particular to reduced flow models and to geophysical flows. A closely related approach is the stochastic advection by Lie transport (SALT) introduced by Holm~\cite{Holm2015} through a variational principle. In SALT, the noise acts on the momentum through a Lie derivative, which contains, besides the advection $(\sigma\cdot\nabla)u$, the stretching term $(\nabla\sigma)^Tu$. In our setting, the momentum is transported by the fast velocity through pure advection, and the limit contains the advection noise $P((\sigma\cdot\nabla)u)$ of the LU model, with $P$ the Leray projection onto divergence-free fields.

A white noise in time is an idealization, since it corresponds to an infinite separation of time scales. In practice, the unresolved velocity has a small but positive correlation time. We believe that it is important to prove that the LU equations are limits of equations with more realistic noises, and this is the aim of this article.

Let $S\subset\R^d$, $d\in\{2,3\}$, be a smooth bounded domain, and let $m(t)$ be a centered, stationary, divergence-free Markov process with values in a space of smooth vector fields on $S$ vanishing on $\partial S$. The fast velocity is of the form $\eps^{-1}m(t/\eps^2)$: its amplitude is of order $\eps^{-1}$ and its correlation time of order $\eps^2$. Along a Lagrangian trajectory, we write
\begin{align}\label{eq:lagrangian-intro}
 \dot X^\eps(t)&=u^\eps(t,X^\eps(t))+\frac1\eps m\!\left(\frac{t}{\eps^2},X^\eps(t)\right)\nonumber\\
 &=u^\eps(t,X^\eps(t))+u_s(X^\eps(t))+\left[\frac1\eps m\!\left(\frac{t}{\eps^2},X^\eps(t)\right)-u_s(X^\eps(t))\right].
\end{align}
The field $u_s$ is the It\^o--Stokes drift. It comes from the spatial inhomogeneity of $m$, and it is subtracted so that the last bracket is asymptotically centered; its expression is given in Section~\ref{sec:fast}. In Eulerian variables, we obtain
\begin{equation}\label{eq:eps-abstract-intro}
 \partial_t u^\eps+B\!\left(u^\eps+\frac1\eps m^\eps,u^\eps+u_s\right)=-Au^\eps-\frac1\eps A m^\eps,\qquad u^\eps(0)=u_0,
\end{equation}
where $m^\eps(t)=m(t/\eps^2)$, $A$ is the Stokes operator with Dirichlet boundary conditions and $B$ is the usual projected nonlinear operator.

Our main result can be stated informally as follows. Precise assumptions are given in Section~\ref{sec:fast}, and the precise statement is Theorem~\ref{thm:main}.

\begin{theorem*}
Assume that $m$ has enough finite moments, that its mean and its two-time correlations decay sufficiently fast, and that its time-integrated correlations are symmetric. Then the solutions $u^\eps$ of \eqref{eq:eps-abstract-intro} obtained by Galerkin approximation are tight as $\eps\to0$, and every limit point is a martingale solution of the LU stochastic Navier--Stokes equation
\begin{equation}\label{eq:limit-intro}
 \dd u+[Au+B(u,u+u_s)+F(u)]\dd t=G(u)\dd W_t,\qquad u(0)=u_0,
\end{equation}
where $F$ and $G$ are given by the time-integrated covariance of $m$ and $W$ is a cylindrical Wiener process. In dimension two, the solution of \eqref{eq:limit-intro} is pathwise unique and $u^\eps$ converges in law to it. Moreover, the convergence holds jointly with the integrated fast velocity $\eps^{-1}\int_0^\cdot m^\eps(s)\dd s$, whose limit is $-Q^{1/2}W$, and the same results hold for Galerkin approximations $u^\eps_{N(\eps)}$ with any $N(\eps)\to\infty$.
\end{theorem*}

Equation \eqref{eq:limit-intro} is the It\^o form of a Stratonovich transport equation. The amplitude $\eps^{-1}$ and the correlation time $\eps^2$ balance each other, so that a transport noise survives in the limit. Its covariance is the Green--Kubo covariance of $m$, and the spatial variations of this covariance produce the It\^o--Stokes drift, which appears in the It\^o correction as an advection by $-u_s$. Thus the drift and the diffusion of the LU model are not chosen independently: both come from the same unresolved velocity. Note that $m$ need not be Gaussian: at leading order, only its integrated covariance matters.

Note that both the large scales, represented by $u$, and the small scales, represented by $m$, converge: the integrated fast velocity $\eps^{-1}\int_0^\cdot m^\eps(s)\dd s$ converges, jointly with $u^\eps$, to $-Q^{1/2}W$. This justifies both the LU modeling and the representation of the small scales by a white noise.

Our basic tool is the perturbed test function method. It was introduced by Papanicolaou, Stroock and Varadhan~\cite{PSV} in finite dimension, and this work is part of a program on diffusion approximation for infinite dimensional systems. Let us mention for instance~\cite{DebusscheVovelle2012} for stochastic kinetic limits, \cite{DebusscheRoselloVovelle2022} for kinetic-fluid systems with mixing Markov forcing, \cite{BarrueBouardDebussche2024} for singular dispersive limits, or~\cite{BarrueDebusscheTusseau2024} for Schr\"odinger equations with a random potential. The terms of size $\eps^{-1}$ in \eqref{eq:eps-abstract-intro} prevent any estimate uniform in $\eps$ from the energy identity, and the generator of the pair $(u^\eps,m^\eps)$ contains terms of order $\eps^{-2}$ and $\eps^{-1}$. Both problems are solved by adding correctors to the test functions; these are obtained from Poisson equations for the generator $M$ of the fast process. We use this idea twice: first on the energy, to get bounds uniform in $\eps$, and then on cylindrical observables, to prove tightness and identify the limit. We do not assume that $m$ is bounded. As in~\cite{BarrueDebusscheTusseau2024}, the correctors are then not uniformly bounded, and we use a stopping time to control the large values of the fast field. Also, a second corrector is already needed for the energy, in order to replace the order-one terms which are quadratic in the fast variable by their averages.

The main difficulty is due to the boundary. The correctors which compensate the transport term are linear in $u$, with coefficients such as $P((n\cdot\nabla)h)$. These fields are divergence free, but the Leray projector does not preserve the Dirichlet condition, so that they do not belong to $D(A^{1/2})$. The perturbed generator then contains terms like $(Au,P((n\cdot\nabla)h))$, which cannot be controlled by the energy, see Section~\ref{sec:formal}. We overcome this problem by using \emph{modified correctors}, obtained by inverting $M-\eps^2A$ instead of $M$. The Stokes operator then acts on the corrector, which satisfies the boundary condition, and the bad terms cancel in It\^o's formula (Corollary~\ref{cor:cancel}). The counterpart is that the modified correctors have a boundary layer of thickness $\eps$, the distance over which viscosity acts during one correlation time $\eps^2$. Their $V$-norm is of order $\eps^{-1+2\theta}$ for every $\theta<1/4$, roughly $\eps^{-1/2}$. This is compensated by the powers of $\eps$ in front of the remainders, so that the boundary layer does not affect the limit, but we need a stronger moment condition on $m$ ($4p<q$ in Assumption~\ref{ass:fast}). On the torus, or whenever the correctors can be taken in $D(A^{1/2})$, the remainders are of order $\eps$. The modified correctors may be compared with the frequency cut-off used in~\cite{DebusscheHofmanova2026Fluct}.

In~\cite{DebusscheHugMemin2023}, we have constructed martingale solutions of the LU Navier--Stokes equation in dimensions two and three, proved pathwise uniqueness in dimension two, and studied the vanishing-noise limit. The limit equation obtained here is this LU equation, up to the additional drift $-\frac12\sum_jB(\sigma_j,A\sigma_j)$. This term is the It\^o--Stratonovich correction of the interaction between the transport noise and the additive noise $AQ^{1/2}\dd W$, which comes from the viscous force acting on the fast component of the velocity; it is not taken into account in the LU model (see Remark~\ref{rem:modelling}). It vanishes for spatially homogeneous noise on the torus (Section~\ref{sec:homogeneous}), so that the difference with~\cite{DebusscheHugMemin2023} is an effect of the inhomogeneity of the noise. We do not know which modeling is the most relevant. We give a self-contained proof of uniqueness in dimension two (Proposition~\ref{prop:LU2D}).

Stochastic Navier--Stokes equations with transport noise go back to Brze\'zniak, Capi\'nski and Flandoli~\cite{BrzezniakCapinskiFlandoli1992}, whose framework covers two-dimensional equations with noise depending on first-order derivatives of the velocity. Mikulevi\v{c}ius and Rozovskii derived stochastic fluid equations from stochastic particle trajectories~\cite{MikuleviciusRozovskii2004} and proved existence and uniqueness for a large class of such equations~\cite{MikuleviciusRozovskii2005}. The LU framework has the same transport structure, and makes the covariance-induced diffusion and the It\^o--Stokes drift explicit. Transport noise has also been obtained from scale separation in two-dimensional fluid models by Flandoli and Pappalettera~\cite{FlandoliPappalettera2021,FlandoliPappalettera2022}, and Debussche and Pappalettera~\cite{DebusschePappalettera} have obtained transport noise and an It\^o--Stokes drift for coupled slow--fast systems, including Navier--Stokes, surface quasi-geostrophic and primitive equations. A rough-path approach to a related two-scale system is given in~\cite{DebusscheHofmanova}. More recently, randomly advected Navier--Stokes equations with both a fast time scale $\eps^2$ and a vanishing correlation length $\delta$ have been studied in~\cite{DebusscheHofmanova2026Fluct}. When $\eps=o(\delta)$, the limit is deterministic, with an enhanced viscosity in which the Leray projector plays a role, and Gaussian fluctuations appear in dimension two. Our setting is complementary since the spatial covariance is fixed and the noise survives at leading order. Navier--Stokes equations with transport noise have also been studied directly, see for instance~\cite{AgrestiVeraar2024,Pappalettera2024,FlandoliLuongo2023}.

It is also known that transport noise has regularizing and dissipative effects. Flandoli and Luo~\cite{FlandoliLuo2021} have shown that high-mode transport noise improves the control of the vorticity in three-dimensional Navier--Stokes equations and gives well-posedness with high probability (see also~\cite{FlandoliGaleatiLuo2021}). Agresti and Veraar~\cite{AgrestiVeraar2024} have developed a theory in critical spaces, with stochastic Serrin-type blow-up criteria. Agresti has shown that transport noise can induce anomalous enstrophy dissipation in two dimensions, even for small noise intensity~\cite{Agresti2025Anomalous}, and that suitable transport noise gives global smooth solutions with high probability for hyperviscous Navier--Stokes equations with dissipation exponent below the Lions threshold $5/4$~\cite{Lions1969,Agresti2024Hyperviscous}, and also for the three-dimensional Navier--Stokes equations~\cite{Agresti2026Absence}. Regularization is not our purpose, but these results show that the effective diffusion created by transport noise is important for the dynamics. In our setting, it is not postulated but obtained in the limit.

Compared with coupled slow--fast systems, the unresolved velocity is here an abstract stationary Markov field, and not the solution of a second fluid equation. We only assume polynomial moments and a polynomially weighted mixing of its mean and of its two-time correlations, and the limiting covariance is written directly in terms of its Green--Kubo tensor. The covariance symmetry \eqref{eq:symcorr} could be relaxed, at the price of an additional drift in the limit. We also think that similar results hold for other geophysical models for which LU versions are available.

The article is organized as follows. In Section~\ref{sec:model}, we derive the model and introduce the functional setting. The assumptions on the fast process and the limiting operators are given in Section~\ref{sec:fast}, and the main result is stated in Section~\ref{sec:main}. Section~\ref{sec:existence} deals with fixed $\eps$. The modified correctors are introduced in Section~\ref{sec:resolvent}. In Section~\ref{sec:energy}, we explain why we work with Galerkin approximations and prove the uniform estimates. In Sections~\ref{sec:ptf} to~\ref{sec:passage}, we construct the correctors for cylindrical observables, identify the limit, prove tightness and conclude. Finally, the spatially homogeneous case is discussed in Section~\ref{sec:homogeneous}.

Throughout the article, $C$ denotes a constant whose value may change from one line to the other and which, unless explicitly stated, does not depend on $\eps$ nor on the Galerkin parameter $N$.

\section{Modeling and functional setting}\label{sec:model}

\subsection{Lagrangian description with finite-correlation random advection}
The starting point is the Lagrangian description underlying modeling under location uncertainty. In the LU framework the trajectory of a fluid particle is decomposed into a smooth-in-time large-scale velocity and a rapidly fluctuating velocity which, in the idealized model, is white in time,
\begin{equation}\label{eq:LU-lagrangian}
 \dd X(t)=v(t,X(t))\,\dd t+\sigma(X(t))\,\dd W_t.
\end{equation}
The white-in-time component represents an infinite separation of temporal scales. At finite scale separation, however, the unresolved velocity retains a non-zero correlation time. We therefore replace the white-noise velocity by a centered stationary random field $m$ with a short correlation time. Under the diffusion scaling the Lagrangian trajectory is described by
\begin{equation}\label{eq:lagrangian}
 \dot X^\eps(t)=u^\eps(t,X^\eps(t))+\frac1\eps m\!\left(\frac{t}{\eps^2},X^\eps(t)\right),\qquad \eps>0.
\end{equation}
The process $m$ is divergence free and sufficiently regular in space; its precise assumptions are given in Section~\ref{sec:fast}. The correlation time is of order $\eps^2$ while the amplitude is of order $\eps^{-1}$.

Although $m(t,x)$ is centered for each fixed $x$, the random field evaluated along the trajectory, $m(t/\eps^2,X^\eps(t))$, is not centered. Since the small scales should average to zero, a correction has to be added. The coupling between the spatial dependence of $m$ and the displacement generated by $m$ itself produces, at order one, the drift which becomes the It\^o--Stokes drift of the limiting LU model. We reproduce the formal argument leading to this correction; its purpose is to identify the physical structure of the model.

\subsection{Formal emergence of the It\^o--Stokes drift}
For clarity, suppose that $m$ has a dependence range in time bounded by one, i.e.\ $(m(r))_{r\le t-1}$ and $(m(r))_{r\ge t}$ are independent for every $t$, and that the initial position is independent of $m$. Since $X^\eps(t-\eps^2)$ depends on $m$ only through $(m(r))_{r\le t/\eps^2-1}$, we have $\E[m(t/\eps^2,X^\eps(t-\eps^2))]=0$ and
\[
 \frac1\eps\E[m(\frac{t}{\eps^2},X^\eps(t))]
 \simeq\frac1\eps\E[\nabla m(\frac{t}{\eps^2},X^\eps(t))\cdot(X^\eps(t)-X^\eps(t-\eps^2))].
\]
In the increment of $X^\eps$, the contribution of the large-scale velocity is of order $\eps$ and vanishes formally. The leading contribution is generated by the fast field itself. Since $X^\eps$ has a time derivative of order $\eps^{-1}$ and $t-r$ is of order $\eps^2$, we may freeze the spatial argument, $m(r/\eps^2,X^\eps(r))=m(r/\eps^2,X^\eps(t))+O(\eps)$. After the change of variable $r=t-\eps^2s$,
\[
 \frac1\eps\E[m_i(\frac{t}{\eps^2},X^\eps(t))]\simeq
 \sum_{j=1}^d\int_0^1\E[\partial_{x_j}m_i(\frac{t}{\eps^2},X^\eps(t))\,m_j(\frac{t}{\eps^2}-s,X^\eps(t))]\dd s .
\]
Formally, we freeze the spatial argument $x=X^\eps(t)$. By stationarity, $\E[\partial_{x_j}m_i(\frac{t}{\eps^2},x)\,m_j(\frac{t}{\eps^2}-s,x)]=\E_\eta[\partial_{x_j}m_i(s,x)\,m_j(0,x)]$ depends only on the time shift $s$ and on $x$, and it vanishes for $s>1$ since $m$ is centered with dependence range one. Using stationarity once more, and the derivative in $x$ of the covariance symmetry \eqref{eq:symcorr} below to reverse the order of the two times, we are led to define the It\^o--Stokes drift by
\begin{equation}\label{eq:us-correlation-gradient}
 (u_s)_i(x):=\sum_{j=1}^d\int_0^\infty\E_\eta[\partial_{x_j}m_i(0,x)\,m_j(s,x)]\dd s
 =\sum_{j=1}^d\partial_{x_j}\int_0^\infty\E_\eta[m_i(0,x)m_j(s,x)]\dd s,
\end{equation}
where the second equality uses $\nabla\cdot m=0$. In terms of the Green--Kubo tensor $k$ of \eqref{eq:k} and the variance tensor $a(x)=k(x,x)$, this reads
\begin{equation}\label{eq:formal-us-covariance}
 u_s=\frac12\nabla\cdot a,\qquad\text{and}\qquad
 \frac1\eps\E[m(\frac{t}{\eps^2},X^\eps(t))]\simeq\E[u_s(X^\eps(t))].
\end{equation}
Since the small scales should have zero average, this quantity is subtracted from the fast part of the velocity, and \eqref{eq:lagrangian} is rewritten as
\begin{equation}\label{eq:lagrangian-us}
 \dot X^\eps(t)=[u^\eps(t,X^\eps(t))+u_s(X^\eps(t))]+[\frac1\eps m(\frac{t}{\eps^2},X^\eps(t))-u_s(X^\eps(t))].
\end{equation}
The first bracket is interpreted as the large-scale velocity, whereas the second one is an asymptotically centered small-scale velocity. In the diffusion limit it is expected to converge to the It\^o transport noise in \eqref{eq:LU-lagrangian}. Note that $u_s$ is in general not divergence free, since $\nabla\cdot u_s=\frac12\sum_{i,j}\partial_i\partial_ja^{ij}$. This is consistent with its role below: $u_s$ enters only through the transported quantity and through the It\^o correction, never as an advecting velocity in the prelimit equation.

\subsection{Transport of a quantity by the random flow}
We now derive the Eulerian equation corresponding to \eqref{eq:lagrangian}, using the Reynolds-transport methodology of LU. Let $V(t)$ be a volume transported by the random flow $X^\eps$, and let $\Phi$ be a test function transported by the same flow, with $\Phi(0)=g$ supported in $V(0)$. We use the ansatz $\partial_t\Phi=\alpha+\eps^{-1}\beta(t/\eps^2,\cdot)$. Since $\Phi$ is conserved along trajectories, \eqref{eq:lagrangian} gives $\alpha=-(u^\eps\cdot\nabla)\Phi$ and $\beta=-(m\cdot\nabla)\Phi$. For a vector-valued quantity $q$, using $\nabla\cdot m^\eps=0$ and integrating by parts,
\begin{equation}\label{eq:transport-random}
 \frac{\dd}{\dd t}\int_{V(t)}q\Phi\,\dd y=\int_S[\partial_tq+(u^\eps+\frac1\eps m^\eps)\cdot\nabla q]\Phi\,\dd y+\int_S(\nabla\cdot u^\eps)q\Phi\,\dd y.
\end{equation}

\subsection{Mass and momentum conservation}
For an incompressible homogeneous fluid with constant density $\rho$, taking $q=\rho$ in \eqref{eq:transport-random} gives $\nabla\cdot u^\eps=0$. We set without loss of generality $\rho=1$. We then apply the transport identity to the large-scale momentum $q=u^\eps+u_s$ identified in \eqref{eq:lagrangian-us}. Newton's second law, with molecular viscosity $\nu$ and a pressure force decomposed into slow and fast components $p+\eps^{-1}p^\eps$, gives
\begin{equation}\label{eq:NS-random-strong}
 \partial_tu^\eps+(u^\eps+\frac1\eps m^\eps)\cdot\nabla(u^\eps+u_s)=\nu\Delta(u^\eps+\frac1\eps m^\eps)-\nabla p-\frac1\eps\nabla p^\eps,
\end{equation}
with $\nabla\cdot u^\eps=\nabla\cdot m^\eps=0$ and $u^\eps=m^\eps=0$ on $\partial S$. Since $u_s$ does not depend on time, $\partial_t(u^\eps+u_s)=\partial_tu^\eps$. The transported quantity is the large-scale momentum, whereas the viscous force acts on the full velocity $u^\eps+\eps^{-1}m^\eps$, which is also the advecting velocity.

\begin{remark}\label{rem:closure}
 The fast field $\eps^{-1}m^\eps$ is prescribed: it is not required to satisfy a momentum balance of its own. In particular, its own acceleration $\partial_t(\eps^{-1}m^\eps)+\eps^{-2}(m^\eps\cdot\nabla)m^\eps$ and the corresponding pressure do not appear in \eqref{eq:NS-random-strong}; they are implicitly balanced by unresolved small-scale forcing, which is not modeled. This is the closure underlying the LU framework~\cite{Memin2014,ResseguierMeminChapron2017}, where the statistics of the unresolved velocity are a modeling input. It differs from coupled slow--fast systems such as~\cite{DebusschePappalettera}, in which the fast component solves its own fluid equation.
\end{remark}

\begin{remark}\label{rem:transported}
The momentum balance is written for $q=u^\eps+u_s$, the large-scale velocity of \eqref{eq:lagrangian-us}. If one transports $q=u^\eps$ instead, $u_s$ disappears from the second argument of $B$ in \eqref{eq:eps}. Theorem~\ref{thm:main} then holds with the same proof, $u+u_s$ being replaced by $u$ in \eqref{eq:Fexplicit}--\eqref{eq:Gexplicit}; the energy estimate even simplifies, since $F_E=0$ in Section~\ref{sec:energy}. In both cases the advection by $-u_s$ in the limit is not inserted by hand: it comes from the It\^o correction, see \eqref{eq:F-decomp}.
\end{remark}

\begin{remark}\label{rem:modelling}
Since the viscous force acts on the full velocity, \eqref{eq:eps} contains the singular, centered forcing $-\eps^{-1}Am^\eps$. In the limit, it produces the additive noise $AQ^{1/2}\dd W$ on the right-hand side of \eqref{eq:LU}, as in the LU model, and the contribution $-\frac12\sum_jB(\sigma_j,A\sigma_j)$ to $F(u)$ in \eqref{eq:Fexplicit}, which appears on the left-hand side of \eqref{eq:LU}. This last term is the It\^o--Stratonovich correction of the interaction between the transport noise and the additive noise; it is not taken into account in the LU model. If the viscous force is applied to the resolved velocity only, these terms disappear: one removes $A\nN$ from the fast drift $\sff$ defined at the end of Section~\ref{sec:existence}, and the proof is unchanged and shorter. We do not know which modeling is the most relevant; all the variants of Remarks~\ref{rem:transported} and~\ref{rem:modelling} are covered by the same arguments.
\end{remark}

\subsection{Divergence-free spaces and the Stokes operator}\label{sec:functional}
Let $\mathcal V$ be the space of smooth compactly supported divergence-free vector fields on $S$. We introduce the closure of this space in two different spaces:
\[
 H=\overline{\mathcal V}^{L^2(S;\R^d)},\qquad V=\overline{\mathcal V}^{H^1_0(S;\R^d)},
\]
and we use the Gelfand triple $V\subset H\subset V'$. We write $|\cdot|$ for the norm of $L^2(S;\R^d)$ (which is the norm of $H$). Let $P$ be the Leray projector and $A=-\nu P\Delta$, $D(A)=H^2(S;\R^d)\cap V$, the Stokes operator~\cite{Temam1983}, and endow $V=D(A^{1/2})$ with the norm $\|u\|=|A^{1/2}u|$, which is equivalent to the $H^1_0$ norm and satisfies $\|u\|^2=(Au,u)$ for $u\in D(A)$. The operator $A$ has an orthonormal eigenbasis $(e_k)_{k\ge1}$ in $H$, with eigenvalues $0<\lambda_1\le\lambda_2\le\cdots$. We set $H_N=\mathrm{span}(e_1,\dots,e_N)$ and denote by $\PN$ the orthogonal projection of $H$ onto $H_N$. The bilinear operator is $B(u,v)=P((u\cdot\nabla)v)$, defined for $u\in H$ and vector fields $v$ that need not be divergence free (this is needed for $v=u_s$). We use the following claims.
\begin{enumerate}[label=(C\arabic*)]
\item\label{F1} For $0\le\beta<1/4$, $D(A^\beta)=H\cap H^{2\beta}(S;\R^d)$ with equivalent norms~\cite{FujitaMorimoto,Giga1985}; no boundary condition is involved in this range.
\item\label{F2} $P$ is bounded on $H^s(S;\R^d)$ for every $s\ge0$, and the semigroup $e^{-\tau A}$ is analytic: $|A^\beta e^{-\tau A}|_{\mathcal L(H)}\le C_\beta\tau^{-\beta}$ and $|(e^{-\tau A}-I)A^{-\beta}|_{\mathcal L(H)}\le C_\beta\tau^\beta$ for $\beta\in[0,1]$, $\tau>0$.
\item\label{F3} $|u|_{L^4}\le C|u|^{1/4}\|u\|^{3/4}$ for $u\in V$ and $d\le3$, and $\|u\|^2\ge\lambda_1|u|^2$.
\item\label{F4} If $w\in H\cap L^3(S;\R^d)$, then $((w\cdot\nabla)z,y)=-((w\cdot\nabla)y,z)$ for $y,z\in H^1(S;\R^d)$, with no boundary term, see~\cite{Temam1983}. In particular, for $u\in V$, $(B(u,v),w)_H=-((u\cdot\nabla)w,v)$ for $v\in H^1$ and $w\in H\cap H^1$, and $(B(u,v),v)_H=0$ if $v\in V$.
\end{enumerate}
We fix $\gamma>5/4$, so that $D(A^{\gamma+1})\subset W^{3,\infty}$, $D(A^\gamma)\subset W^{1,\infty}$, and $B$ extends continuously from $H\times L^2$ to $D(A^{-\gamma})$:
\begin{equation}\label{eq:Bnegative}
 \|B(u,v)\|_{D(A^{-\gamma})}\le C|u||v|.
\end{equation}
After applying the Leray projector, \eqref{eq:NS-random-strong} becomes
\begin{equation}\label{eq:eps}
 \partial_tu^\eps+B(u^\eps+\frac1\eps m^\eps,u^\eps+u_s)=-Au^\eps-\frac1\eps Am^\eps,\qquad u^\eps(0)=u_0\in H.
\end{equation}
Our aim is to investigate the limit $\eps\to0$ of this equation.

\begin{remark}\label{rem:gamma}
The condition $\gamma>5/4$, which will lead to the very smooth state space $E=D(A^{\gamma+1})$ for the fast field in Section~\ref{sec:fast}, is certainly not optimal. It is used to ensure that the covariance $k$ belongs to $C^3$, that the term $A^2\PN L_1$ appearing in the energy estimate is bounded in $H$, that $D(A^{\gamma+1/2})\subset W^{2,\infty}$ in the estimate of the Galerkin error, and \eqref{eq:Bnegative}. We have not tried to optimize it. This is not important for our purpose, since in the LU modeling the noise is always very regular in space.
\end{remark}

\section{The fast Markov advection and its effective covariance}\label{sec:fast}

\subsection{Assumptions}

Let $E:=D(A^{\gamma+1})$, a Hilbert space for the norm $\|n\|_E=|A^{\gamma+1}n|$; its elements are divergence free and vanish on $\partial S$. For $p\ge0$ we set $\mathrm w_p(n):=1+\|n\|_E^p$. We denote by $\mathcal L_1(E)$ the space of trace-class operators on $E$, with norm $\|\cdot\|_{\mathcal L_1(E)}$, and, for $x,y\in E$, by $x\otimes y$ the rank-one operator $z\mapsto(y,z)_Ex$. Recall that every $K\in\mathcal L_1(E)$ has a singular value decomposition $K=\sum_is_i\,x_i\otimes y_i$, where $(x_i)$ and $(y_i)$ are orthonormal families in $E$ and $s_i\ge0$, with $\sum_is_i=\|K\|_{\mathcal L_1(E)}$. Hence every bounded bilinear map $\mathcal B:E\times E\to X$ with values in a Banach space $X$ extends uniquely to a bounded linear map on $\mathcal L_1(E)$, with the same norm, by setting $\mathcal B(K)=\sum_is_i\,\mathcal B(x_i,y_i)$.

\begin{assumption}\label{ass:fast}
\begin{enumerate}[label=(\alph*)]
\item \emph{(Markov property and moments.)} $m=(m(t))_{t\in\R}$ is a centered, stationary Markov process with respect to a filtration $(\mathcal F_t)$, with values in a measurable subset $\mathscr E\subset E$, its state space, with c\`adl\`ag paths in $E$, a measurable transition function, transition semigroup $(\mathsf P_t)_{t\ge0}$, generator $M$ and invariant law $\eta$. We write $\widetilde m(\cdot,n)$ for the process started at $n\in\mathscr E$. For some $q>0$ (constrained by \eqref{eq:exponents} below),
\begin{equation}\label{eq:moments}
 \E[\sup_{0\le t\le1}\|m(t)\|_E^q]<\infty .
\end{equation}
\item \emph{(First-order mixing.)} There exist $r_1\ge0$ and a bounded function $\varrho_0\ge0$ with $\int_0^\infty(1+t)\varrho_0(t)\dd t<\infty$ such that
\begin{equation}\label{eq:mix1}
 \|\E\,\widetilde m(t,n)\|_E\le\varrho_0(t)\,\mathrm w_{r_1}(n),\qquad t\ge0,\ n\in\mathscr E.
\end{equation}
\item \emph{(Second-order mixing.)} There exist $r_2\ge0$ and bounded functions $\varrho_1,\varrho_2\ge0$ with $\int_0^\infty(1+t)(\varrho_1+\varrho_2)(t)\dd t<\infty$ such that, for all $s,t\ge0$ and $n\in\mathscr E$,
\begin{equation}\label{eq:mix2}
 \Bigl\|\E[\widetilde m(t,n)\otimes\widetilde m(t+s,n)]-\E_\eta[m(0)\otimes m(s)]\Bigr\|_{\mathcal L_1(E)}\le\varrho_1(t)\varrho_2(s)\,\mathrm w_{r_2}(n).
\end{equation}
\item \emph{(Compatibility of exponents.)} With $p:=1+\max(r_1,r_2)$,
\begin{equation}\label{eq:exponents}
 4p<q .
\end{equation}
\item The covariance symmetry \eqref{eq:symcorr} holds.
\end{enumerate}
\end{assumption}

The process $m$ is \emph{not} assumed to be bounded. The estimates are only required on the state space $\mathscr E$, which carries $\eta$ and the paths of $m$; all the functions of $n$ considered below are defined on $\mathscr E$. Note that, at $t=0$, \eqref{eq:mix1} gives $\|n\|_E\le\varrho_0(0)\mathrm w_{r_1}(n)$ on $\mathscr E$, so that $r_1\ge1$ unless $\mathscr E$ is bounded. Condition \eqref{eq:mix1} quantifies the decay of the mean, and \eqref{eq:mix2} the decay of the two-time correlations, starting from a deterministic state; the polynomial weights allow unbounded state spaces. By \eqref{eq:exponents} we may fix once and for all
\begin{equation}\label{eq:theta-alpha}
 \theta\in(\frac pq,\frac14),\qquad\alpha\in(\frac2q,\frac{2\theta}p),
\end{equation}
so that \ref{F1} applies with $\beta=\theta$. All the moments used below are of order at most $4p<q$.

By \eqref{eq:mix1}, the first corrector of the fast process,
\begin{equation}\label{eq:L1}
 L_1(n):=-\int_0^\infty\E\,\widetilde m(t,n)\dd t,\qquad\text{so that }ML_1=\mathrm{Id},
\end{equation}
is well defined, divergence free, and $\|L_1(n)\|_E\le\|\varrho_0\|_{L^1}\mathrm w_{r_1}(n)$.

\begin{lemma}[Large values of the fast process]\label{lem:stop}
For $T>0$ and $\eps\in(0,1)$, set $m^\eps(t)=m(t/\eps^2)$ and
\begin{equation}\label{eq:tau}
 \tau^\eps:=\inf\bigl\{t\in[0,T]:\|m^\eps(t)\|_E\ge\eps^{-\alpha}\bigr\},\qquad\inf\emptyset=T .
\end{equation}
Then $\|m^\eps(t)\|_E<\eps^{-\alpha}$ for $t<\tau^\eps$, $\Pp(\tau^\eps<T)\le C(T+1)\eps^{\alpha q-2}\to0$, and for $p'\le q$,
\begin{equation}\label{eq:supmoment}
 \E\sup_{t\le T}\|m^\eps(t)\|_E^{p'}\le C(T+1)^{p'/q}\eps^{-2p'/q} .
\end{equation}
\end{lemma}

\begin{proof}
Let us cover $[0,T\eps^{-2}]$ by $K_\eps\le T\eps^{-2}+1$ unit intervals and let $Y_k$ be the supremum of $\|m\|_E$ on the $k$-th one. By stationarity, the $Y_k$ have the law of $Y_0=\sup_{[0,1]}\|m\|_E$.

The union bound and Markov's inequality give
\[
 \Pp\big(\max_kY_k\ge\eps^{-\alpha}\big)\le K_\eps\eps^{\alpha q}\E Y_0^q\le C(T+1)\eps^{\alpha q-2},
\]
which tends to $0$ since $\alpha q>2$.

For the moments, we use Jensen's inequality:
\[
 \E\max_kY_k^{p'}\le\Big(\E\max_kY_k^{q}\Big)^{p'/q}\le\Big(\sum_k\E Y_k^q\Big)^{p'/q}=\big(K_\eps\E Y_0^q\big)^{p'/q},
\]
and \eqref{eq:supmoment} follows. See also~\cite[Lemma~2.1]{BarrueDebusscheTusseau2024}.
\end{proof}

\begin{example}[Ornstein--Uhlenbeck fields]\label{ex:OU}
Let $Q_0$ be a non-negative operator on $E$, self-adjoint and trace class for the scalar product $(x,y)_E=(A^{\gamma+1}x,A^{\gamma+1}y)$ of $E$, and let $m$ be the stationary solution of $\dd m=-m\dd t+Q_0^{1/2}\dd\beta_t$, where $\beta$ is a cylindrical Wiener process on $E$. Then $m$ is Gaussian with invariant law $\mathcal N(0,\frac12Q_0)$ and continuous paths in $E$, and we take $\mathscr E=E$ and \eqref{eq:moments} holds for every $q$. Since $\E\,\widetilde m(t,n)=e^{-t}n$, \eqref{eq:mix1} holds with $\varrho_0(t)=e^{-t}$ and $r_1=1$, and $L_1(n)=-n$. Moreover
\[
 \E[\widetilde m(t,n)\otimes\widetilde m(t+s,n)]-\E_\eta[m(0)\otimes m(s)]=e^{-s-2t}(n\otimes n-\tfrac12Q_0),
\]
so \eqref{eq:mix2} holds with $\varrho_1(t)=e^{-2t}$, $\varrho_2(s)=e^{-s}$ and $r_2=2$. Hence $p=3$, and \eqref{eq:exponents} holds. The process is reversible, so \eqref{eq:symcorr} holds. Since $\E[m_i(0,x)m_j(t,y)]=e^{-t}\E_\eta[m_i(x)m_j(y)]$, the Green--Kubo kernel is $k^{ij}(x,y)=2\E_\eta[m_i(x)m_j(y)]$, which is continuous because $E\subset C^3(\overline S)$. The corresponding operator $Q$ on $H$ is related to $Q_0$ as follows. Let $J:E\to H$ be the inclusion and $J^*=A^{-2(\gamma+1)}$ its adjoint with respect to the scalar products of $H$ and $E$. Then
\[
 (Qf,g)=2\E_\eta[(m,f)(m,g)]=2\E_\eta[(m,J^*f)_E(m,J^*g)_E]=(Q_0J^*f,J^*g)_E,
\]
that is, $Q=JQ_0J^*$.
In particular $\operatorname{Tr}(AQA)=\|AJQ_0^{1/2}\|^2_{L_2(E,H)}\le\|AJ\|^2_{\mathcal L(E,H)}\operatorname{Tr}_EQ_0$, in accordance with Lemma~\ref{lem:cov}. This covers the spatially smooth Ornstein--Uhlenbeck fields of~\cite{DebusscheHofmanova2026Fluct}, at fixed correlation length, in a bounded domain.

Non-Gaussian examples can be obtained by adapting~\cite[Section~2.3]{BarrueDebusscheTusseau2024}. Let $X(t,x)$ be the solution of a dissipative equation $\dd X=(-A^rX+\mathcal G(X))\dd t+\Sigma\dd\beta_t$ in $H$, with $X(0)=x$. Here $r>d/2$, so that $e^{-tA^r}$ is Hilbert--Schmidt with $\|e^{-tA^r}\|_{L_2(H)}\le Ct^{-d/(4r)}$, $\mathcal G$ is bounded and Lipschitz, $\Sigma$ is invertible on $H$, and $\beta$ is a cylindrical Wiener process. Let $\bar x$ be the mean of its invariant measure, and set $\widetilde m(t,n)=A^{-\gamma-1}\bigl(X(t,A^{\gamma+1}n+\bar x)-\bar x\bigr)$ for $n\in E$. Since $A^{\gamma+1}$ is an isometry from $E$ onto $H$, this defines a Markov process in $E$, with $\mathscr E=E$, whose values are divergence free and vanish on $\partial S$. As shown in~\cite[Section~2.3 and Appendix]{BarrueDebusscheTusseau2024}, it has polynomial moments of all orders, and for $n_1,n_2\in E$ there is a coupling $(m_1,m_2)$ of $\widetilde m(\cdot,n_1)$ and $\widetilde m(\cdot,n_2)$ with
\[
 \Pp\bigl(m_1(t)\ne m_2(t)\bigr)\le C\bigl(1+\|n_1\|_E+\|n_2\|_E\bigr)^{\kappa_0}e^{-ct}.
\]
This coupling estimate implies \eqref{eq:mix1}--\eqref{eq:mix2}. Indeed, if $f$ is centered with $|f|\le C\mathrm w_{p'}$, the Cauchy--Schwarz inequality gives
\[
 |\mathsf P_tf(n)|\le\int_E\E\bigl[|f(m_1(t))-f(m_2(t))|\bigr]\eta(\dd n_2)\le C\,\mathrm w_{p''}(n)\,e^{-ct/2}
\]
for some $p''$. Applying this to $f(n)=n$ gives \eqref{eq:mix1}, and applying it to $\phi_s(n)=n\otimes\E\,\widetilde m(s,n)$, whose size is bounded by $C\mathrm w_{1+r_1}(n)e^{-cs/2}$, gives \eqref{eq:mix2} with $\varrho_1(t)=e^{-ct/2}$ and $\varrho_2(s)=e^{-cs/2}$. Since all moments are finite, \eqref{eq:exponents} holds as soon as \eqref{eq:moments} is checked with $q$ large, which follows from the same energy estimates. Thus this example satisfies (a)--(d) of Assumption~\ref{ass:fast}. Condition (e) is an additional assumption. It holds for instance if the equation for $X$ is of gradient type, that is $\Sigma=I$ and $\mathcal G=-\nabla G_0$ for some potential $G_0$: then $X$, and hence $m$, is reversible, see~\cite{DaPratoZabczyk1996}.
\end{example}

\begin{example}[Finite-state modulation]\label{ex:chain}
Let $w_1,\dots,w_K\in D(A^s)$ for all $s$ be linearly independent, and let $(\xi(t))_{t\in\R}$ be a stationary, irreducible, reversible Markov chain on a finite set $\Xi\subset\R^K$ whose invariant law $\pi$ satisfies $\sum_\zeta\zeta\pi(\zeta)=0$. Set $m(t)=\sum_k\xi_k(t)w_k$. The state space $\mathscr E=\{\sum_k\zeta_kw_k:\zeta\in\Xi\}$ is finite, so that $m$ is bounded, and \eqref{eq:moments} holds for every $q$. Irreducibility gives $\sup_n\|\mathsf P_t(n,\cdot)-\eta\|_{TV}\le Ce^{-\lambda_*t}$, hence \eqref{eq:mix1} with $r_1=0$. Writing $\E[\widetilde m(t)\otimes\widetilde m(t+s)]=\mathsf P_t\phi_s(n)$ with $\phi_s(n')=n'\otimes\E\,\widetilde m(s,n')$, we get \eqref{eq:mix2} with $\varrho_1=\varrho_2=Ce^{-\lambda_*\cdot}$ and $r_2=0$. Reversibility gives \eqref{eq:symcorr}, and $Q$ has finite rank.
\end{example}

\subsection{Green--Kubo covariance, It\^o--Stokes drift and limiting operators}
For $x,y\in S$ define
\begin{equation}\label{eq:k}
 k^{ij}(x,y)=\int_0^\infty\E[m_i(0,x)m_j(t,y)+m_i(t,x)m_j(0,y)]\dd t=\int_{\R}\E[m_i(0,x)m_j(t,y)]\dd t,
\end{equation}
the integrals converging by \eqref{eq:mix1}. Thus $k$ is the covariance of the Brownian limit of $\eps^{-1}m(t/\eps^2)$. We impose the reversibility-type symmetry
\begin{equation}\label{eq:symcorr}
 \int_0^\infty\E[m_i(0,x)m_j(t,y)]\dd t=\int_0^\infty\E[m_i(t,x)m_j(0,y)]\dd t,
\end{equation}
then $k^{ij}(x,y)=2\int_0^\infty\E[m_i(0,x)m_j(t,y)]\dd t$. By \eqref{eq:L1} and the Markov property,
\begin{equation}\label{eq:k-L1}
 k^{ij}(x,y)=-2\int_En_i(x)[L_1(n)]_j(y)\,\eta(\dd n),\qquad\text{i.e.}\qquad\E_\eta[L_1(n)\otimes n]=-\tfrac12k
\end{equation}
as elements of $\mathcal L_1(E)$, which is meaningful since $\E_\eta[\|n\|_E\mathrm w_{r_1}(n)]<\infty$. Here an element $K=\sum_is_i\,x_i\otimes y_i$ of $\mathcal L_1(E)$ is identified with the continuous kernel $(x,y)\mapsto\sum_is_i\,x_i(x)\otimes y_i(y)$. This identification is injective: if the kernel vanishes, then $(KJ^*Pg,J^*Pf)_E=0$ for all $f,g\in L^2(S;\R^d)$, where $J^*=A^{-2(\gamma+1)}$ is the adjoint of the embedding $J:E\to H$, and $J^*P(L^2)$ is dense in $E$. Let $Q:H\to H$, $Qf(y)=\int_Sk(x,y)f(x)\dd x$. The kernel $k$ is symmetric and non-negative: for $f\in H$, stationarity and the integrability of the correlations give
\[
 (Qf,f)=\lim_{T\to\infty}\frac1T\E\Bigl|\int_0^T(m(t),f)\dd t\Bigr|^2\ge0 .
\]
Let $(\mu_j)_{j\ge0}$ be the positive eigenvalues of $Q$ and $(\varphi_j)_{j\ge0}$ an orthonormal family of corresponding eigenvectors, so that $Q=\sum_j\mu_j(\varphi_j,\cdot)\varphi_j$. We set $\sigma_j=Q^{1/2}\varphi_j=\sqrt{\mu_j}\varphi_j$. The kernel of $Q$ may be non-trivial (for instance, $Q$ has finite rank in Example~\ref{ex:chain}); it plays no role, since $Q^{1/2}$ vanishes on it.

\begin{lemma}[Regularity of the covariance]\label{lem:cov}
Under Assumption~\ref{ass:fast}, $\varphi_j\in E$ for every $j$, and
\begin{equation}\label{eq:noise-smooth}
 \sum_{j\ge0}\mu_j|A\varphi_j|^2<\infty,\qquad\sum_{j\ge0}\mu_j\|\varphi_j\|_{L^\infty(S)}^2<\infty .
\end{equation}
Moreover $k\in C^3(\overline S\times\overline S)$, the variance tensor $a$ defined below belongs to $W^{3,\infty}(S)$, and $u_s\in W^{2,\infty}(S;\R^d)$.
\end{lemma}

\begin{proof}
By \eqref{eq:k-L1}, $k=-2\E_\eta[L_1(n)\otimes n]$ is a Bochner integral in $\mathcal L_1(E)$, and
\[
 \|k\|_{\mathcal L_1(E)}\le2\E_\eta\bigl[\|L_1(n)\|_E\|n\|_E\bigr]\le2\|\varrho_0\|_{L^1}\E_\eta\bigl[\|n\|_E\mathrm w_{r_1}(n)\bigr]<\infty,
\]
by \eqref{eq:L1} and since $1+r_1<q$.

Therefore $k$ is a trace-class operator on $E$, and it has a singular value decomposition
\[
 k=\sum_is_i\,x_i\otimes y_i,
\]
where $(x_i)$ and $(y_i)$ are orthonormal families in $E$ and the $s_i\ge0$ are the singular values of $k$. The trace norm is the sum of the singular values, so that $\sum_is_i=\|k\|_{\mathcal L_1(E)}<\infty$ by the above bound.

We deduce that $Qf=\sum_is_i(x_i,f)y_i$ takes values in $E$ for $f\in H$. In particular $\varphi_j=\mu_j^{-1}Q\varphi_j\in E$, since $\mu_j>0$.

Also, $AQA$ extends to a trace-class operator on $H$ with trace norm at most
\[
 \sum_is_i|Ax_i||Ay_i|\le C\sum_is_i,
\]
since $E\subset D(A)$ continuously. Since $Q\ge0$, we obtain
\[
 \sum_j\mu_j|A\varphi_j|^2=\|AQ^{1/2}\|^2_{L_2(H,H)}=\operatorname{Tr}(AQA)<\infty .
\]
The second condition in \eqref{eq:noise-smooth} follows from $\|\varphi\|_{L^\infty}\le C|A\varphi|$, which holds since $H^2\subset L^\infty$ for $d\le3$.

Finally, $E\subset W^{3,\infty}$, so the series for $k$ converges in $C^3(\overline S\times\overline S)$, and $a(x)=k(x,x)$, $u_s=\frac12\nabla\cdot a$ have the stated regularity.
\end{proof}

Since $|A\varphi|\ge\lambda_1|\varphi|$, \eqref{eq:noise-smooth} gives $\sum_j\mu_j<\infty$, so $Q$ is trace class. The variance tensor and the It\^o--Stokes drift are
\begin{equation}\label{eq:a-us}
 a(x)=k(x,x)=\sum_j\mu_j\varphi_j(x)\otimes\varphi_j(x),\qquad u_s=\frac12\nabla\cdot a,
\end{equation}
the series converging uniformly by \eqref{eq:noise-smooth}. The limiting drift and diffusion operators are
\begin{align}
 F(u)&=-\frac12\sum_{j\ge0}[B(\sigma_j,A\sigma_j)+B(\sigma_j,B(\sigma_j,u+u_s))],\label{eq:Fexplicit}\\
 G(u)w&=AQ^{1/2}w+B(Q^{1/2}w,u+u_s),\qquad w\in H.\label{eq:Gexplicit}
\end{align}
By \ref{F4} and \eqref{eq:noise-smooth}, both series in $F$ converge in $V'$: for $h\in V$ and $v\in H^1$,
\begin{gather*}
 |(B(\sigma_j,A\sigma_j),h)|\le\mu_j\|\varphi_j\|_{L^\infty}|A\varphi_j|\|h\|,\qquad
 |(B(\sigma_j,B(\sigma_j,v)),h)|\le\mu_j\|\varphi_j\|_{L^\infty}^2\|v\|_{H^1}\|h\|.
\end{gather*}
Note that the Leray projector inside $B(\sigma_j,B(\sigma_j,v))$ cannot be removed. Since each $\varphi_j$ is divergence free, $\sum_j\mu_j(\varphi_j\cdot\nabla)(\varphi_j\cdot\nabla)v=\nabla\cdot(a\nabla v)$, and therefore
\begin{equation}\label{eq:F-decomp}
 \sum_{j\ge0}B(\sigma_j,B(\sigma_j,v))=P\nabla\cdot(a\nabla v)-\mathcal R(v),\qquad
 \mathcal R(v)=\sum_{j\ge0}\mu_jP[(\varphi_j\cdot\nabla)(I-P)((\varphi_j\cdot\nabla)v)].
\end{equation}
$\mathcal R$ is a non-local pressure correction; compare with the role of the Leray projection in the enhanced viscosity of~\cite{DebusscheHofmanova2026Fluct}. Besides, $\frac12\nabla\cdot(a\nabla v)=\frac12a:\nabla^2v+(u_s\cdot\nabla)v$, so the It\^o correction contains the advection of $v=u+u_s$ by $-u_s$, as in the LU model. Finally, equation \eqref{eq:LU} below is the It\^o form of
\begin{equation}\label{eq:LU-strat}
 \dd u+[Au+B(u,u+u_s)]\dd t+G(u)\circ\dd W_t=0,
\end{equation}
since the It\^o--Stratonovich correction of $-G(u)\circ\dd W$ is $\frac12\sum_jB(\sigma_j,G(u)\varphi_j)=-F(u)$ (see~\cite{FlandoliLuongo2023}); the sign in front of $W$ is irrelevant in law.

\section{Limiting LU equation and main result}\label{sec:main}

Let us now define the notions of solution that we use and state our main result.

Let $W$ be a cylindrical Wiener process on $H$. The limiting equation, in It\^o form, is
\begin{equation}\label{eq:LU}
 \dd u+[Au+B(u,u+u_s)+F(u)]\dd t=G(u)\dd W_t,\qquad u(0)=u_0 .
\end{equation}

\begin{definition}[Probabilistically weak solution of the random-advection equation]\label{def:eps}
For fixed $\eps>0$, a probabilistically weak solution of \eqref{eq:eps} consists of a stochastic basis carrying processes $(u^\eps,m^\eps)$ such that $m^\eps$ has the law of $(m(t/\eps^2))_{t\ge0}$, almost every path of $u^\eps$ lies in $L^\infty(0,T;H)\cap L^2(0,T;V)\cap C([0,T];D(A^{-\gamma}))$, and \eqref{eq:eps} holds in the weak sense.
\end{definition}

\begin{definition}[Martingale solution of the LU equation]
A martingale solution of \eqref{eq:LU} consists of a stochastic basis, a cylindrical Wiener process $W$, and a progressively measurable process $u\in L^2(\Omega\times(0,T);V)\cap L^2(\Omega;L^\infty(0,T;H))$ satisfying \eqref{eq:LU} weakly.
\end{definition}

Let $\varsigma>0$. We set
\[
 \mathcal X_T=C([0,T];D(A^{-\varsigma}))\cap L^2(0,T;H),\qquad\mathcal Z_T=C([0,T];D(A^{-\varsigma})),
\]
and we denote by
\[
 z^\eps(t):=\frac1\eps\int_0^tm^\eps(s)\dd s,\qquad t\in[0,T],
\]
the integrated fast velocity. A family of laws indexed by $\eps$ is said to be tight as $\eps\to0$ if, from every sequence $\eps_k\to0$, one can extract a subsequence along which these laws converge.

\begin{theorem}[Diffusion approximation]\label{thm:main}
Let $d\in\{2,3\}$, let $S$ be a smooth bounded domain, $u_0\in H$, and let Assumption~\ref{ass:fast} hold. Then:
\begin{enumerate}[label=(\roman*)]
\item for every $\eps>0$, equation \eqref{eq:eps} admits a probabilistically weak solution $u^\eps$, obtained as a limit in law of Galerkin approximations; if $d=2$, this solution is pathwise unique, adapted to the filtration of $m^\eps$, and it is the almost sure limit of the whole Galerkin sequence (Proposition~\ref{prop:2D});
\item let $N(\eps)\to\infty$ as $\eps\to0$; then the family $\{\Law(u^\eps_{N(\eps)},z^\eps)\}$ is tight as $\eps\to0$ in $\mathcal X_T\times\mathcal Z_T$, and every limit point is the law of a pair $(u,Z)$, where $Z=-Q^{1/2}W$ for a cylindrical Wiener process $W$ on $H$, and $u$ is a martingale solution of \eqref{eq:LU} driven by $W$;
\item the same holds for the family $\{\Law(u^\eps,z^\eps)\}$;
\item if $d=2$, \eqref{eq:LU} has pathwise unique solutions (Proposition~\ref{prop:LU2D}); consequently $(u^\eps,z^\eps)$ and $(u^\eps_{N(\eps)},z^\eps)$ converge in law in $\mathcal X_T\times\mathcal Z_T$, as $\eps\to0$, to $(u,Z)$, where $Z=-Q^{1/2}W$ and $u$ is the unique solution of \eqref{eq:LU} driven by $W$.
\end{enumerate}
\end{theorem}

Since $G(u)w=AQ^{1/2}w+B(Q^{1/2}w,u+u_s)$, the limit equation can be written
\[
 \dd u+[Au+B(u,u+u_s)+F(u)]\dd t+A\dd Z+B(\dd Z,u+u_s)=0,
\]
which is the It\^o form of $\dd u+[Au+B(u,u+u_s)]\dd t+A\circ\dd Z+B(\circ\dd Z,u+u_s)=0$. Thus the noise of the limit equation is the limit of the integrated fast velocity, as suggested by \eqref{eq:eps}. In (ii), $N(\eps)$ is arbitrary, so that the limits $\eps\to0$ and $N\to\infty$ may be taken simultaneously.

\begin{remark}\label{rem:scope}
 In dimension three, \eqref{eq:eps} is not known to have unique solutions, and Theorem~\ref{thm:main} concerns a selection: the probabilistically weak solutions obtained as limits of Galerkin approximations in Section~\ref{sec:existence}. All the estimates are proved on these approximations, and we do not claim anything about other possible solutions of \eqref{eq:eps}. In dimension two this restriction disappears: for fixed $\eps$, \eqref{eq:eps} is a two-dimensional Navier--Stokes equation with random but smooth coefficients, its solution is pathwise unique, and the theorem applies to it. Tightness is only asserted as $\eps\to0$, because the energy bounds are uniform up to the stopping time $\tau^\eps$ of Lemma~\ref{lem:stop}, and $\Pp(\tau^\eps<T)$ is small only for small $\eps$. No uniqueness of the limit problem is needed for (ii)--(iii); in dimension two it holds, which gives (iv).
\end{remark}

\section{Existence for fixed correlation scale}\label{sec:existence}

We first construct solutions for a fixed $\eps$. This is classical, but the estimates obtained in this way blow up when $\eps\to0$. Uniform estimates are proved in Section~\ref{sec:energy}. The Galerkin equation is
\begin{equation}\label{eq:galerkin}
 \partial_tu_N^\eps+\PN B(u_N^\eps+\frac1\eps\PN m^\eps,u_N^\eps+u_s)=-Au_N^\eps-\frac1\eps A\PN m^\eps,\qquad u_N^\eps(0)=\PN u_0 .
\end{equation}
For fixed $\eps$ and $N$, this is a finite-dimensional differential equation driven by the c\`adl\`ag process $\PN m^\eps$, which has a unique local solution, adapted to the filtration of $m^\eps$. Testing against $u_N^\eps$ and using \ref{F4},
\[
 \frac12\frac{\dd}{\dd t}|u_N^\eps|^2+\|u_N^\eps\|^2=-\frac1\eps(Am^\eps,u_N^\eps)-\frac1\eps(B(\PN m^\eps,u_s),u_N^\eps)-(B(u_N^\eps,u_s),u_N^\eps).
\]
Hence $\frac{\dd}{\dd t}|u_N^\eps|^2+\|u_N^\eps\|^2\le C(1+\eps^{-2}\|m^\eps(t)\|_E^2)(1+|u_N^\eps|^2)$, and Gronwall's lemma gives global existence and the pathwise bound
\[
 \sup_{t\le T}|u_N^\eps(t)|^2+\int_0^T\|u_N^\eps\|^2\dd t\le C_{\eps,T}(1+|u_0|^2),\qquad C_{\eps,T}=C\exp(CT(1+\eps^{-2}\sup_{t\le T}\|m^\eps(t)\|_E^2)),
\]
which is finite almost surely, independent of $N$, but not integrable with respect to $\omega$ in general. On the events $\Omega_{\eps,R}=\{\sup_{t\le T}\|m^\eps(t)\|_E\le R\}$, whose probability tends to $1$ as $R\to\infty$, this bound is deterministic. From \eqref{eq:galerkin} and \eqref{eq:Bnegative}, on $\Omega_{\eps,R}$,
\[
 \|u_N^\eps(t)-u_N^\eps(s)\|_{D(A^{-\gamma})}\le C\int_s^t(\|u_N^\eps\|+|u_N^\eps|^2+1)\dd r+\frac{C_R}\eps\int_s^t(1+|u_N^\eps|)\dd r .
\]
With the compact embedding $V\Subset H$, this gives tightness of the laws of $u_N^\eps$, $N\ge1$, in $L^2(0,T;H)\cap C([0,T];D(A^{-\gamma}))$. Since the law of $m^\eps$ does not depend on $N$, the laws of the pairs $(u_N^\eps,m^\eps)$, $N\ge1$, are tight in the product of this space with the Skorokhod space $D([0,T];E)$, which is Polish. On a Skorokhod representation space for $(u_N^\eps,m^\eps)$, strong convergence in $L^2(0,T;H)$ and weak convergence in $L^2(0,T;V)$ allow passing to the limit in $B(u_N^\eps,u_N^\eps+u_s)$, and the terms containing $\PN m^\eps$ converge pathwise. The pathwise bound passes to the limit by lower semicontinuity. This gives a probabilistically weak solution of \eqref{eq:eps}, as a limit in law of a subsequence $u^\eps_{N_j}$.

\begin{proposition}[Pathwise uniqueness in dimension two]\label{prop:2D}
Let $d=2$ and $\eps>0$. For every c\`adl\`ag path $m^\eps:[0,T]\to E$, equation \eqref{eq:eps} has at most one solution in $L^\infty(0,T;H)\cap L^2(0,T;V)$. Consequently the whole Galerkin sequence $u^\eps_N$ converges, for almost every $\omega$, in $\mathcal X_T$ (with $\varsigma\ge\gamma$) to the unique solution $u^\eps$, which is adapted to the filtration of $m^\eps$.
\end{proposition}

\begin{proof}
Let us fix a path of $m^\eps$. It is bounded in $E\subset W^{3,\infty}$ on $[0,T]$. Let $u$ be a solution in $L^\infty(0,T;H)\cap L^2(0,T;V)$. In dimension two, $|(B(u,u),w)|\le C|u|\|u\|\|w\|$ by Ladyzhenskaya's inequality, so $\partial_tu\in L^2(0,T;V')$. Hence $u\in C([0,T];H)$, and the energy equality holds by the Lions--Magenes lemma~\cite{Temam1983}.

Let $u_1,u_2$ be two solutions and $w=u_1-u_2$. Since $(B(u_2+\eps^{-1}m^\eps,w),w)=0$ by \ref{F4},
\[
 \frac12\frac{\dd}{\dd t}|w|^2+\|w\|^2=-(B(w,u_1+u_s),w)=(B(w,w),u_1)-(B(w,u_s),w).
\]
By Ladyzhenskaya's inequality, $|(B(w,w),u_1)|\le C|w|\|w\|\|u_1\|\le\frac12\|w\|^2+C\|u_1\|^2|w|^2$, and $|(B(w,u_s),w)|\le\|\nabla u_s\|_{L^\infty}|w|^2$. Since $\|u_1\|^2\in L^1(0,T)$, Gronwall's lemma gives $w=0$.

For the Galerkin approximations, the pathwise bound and the time-increment estimate above are deterministic once the path of $m^\eps$ is fixed; hence every subsequence of $(u^\eps_N)$ has a further subsequence converging in $\mathcal X_T$ to a solution of \eqref{eq:eps}, which is unique. The whole sequence therefore converges for every path, and the limit is adapted, as a limit of adapted processes.
\end{proof}

\emph{Notation.} From now on and until the end of Section~\ref{sec:tight}, all the computations are performed on the Galerkin approximations. We fix $N\ge1$, we denote by $u=u_N^\eps$ the solution of the Galerkin equation \eqref{eq:galerkin}, which takes values in $H_N$, and we set $m^\eps(t)=m(t/\eps^2)$. We also write
\[
 \nN=\PN n,\qquad v=u+u_s,\qquad\sff=A\nN+\PN B(\nN,v),\qquad\sfs=Au+\PN B(u,v),
\]
so that $\partial_tu=-\eps^{-1}\sff(u,m^\eps)-\sfs(u)$. The formal generator of $(u,m^\eps)$ is
\begin{equation}\label{eq:generator}
 \Lc^\eps\Phi=\eps^{-2}M\Phi-\eps^{-1}D_u\Phi[\sff]-D_u\Phi[\sfs].
\end{equation}
Test vectors $h$ are taken in $H_{N_0}$ for some fixed $N_0$, and we always take $N\ge N_0$, so that $\PN h=h$. Note that $\bigcup_{N_0}H_{N_0}$ is dense in $D(A^s)$ for every $s$.

\section{Modified correctors}\label{sec:resolvent}

\subsection{Formal construction of the correctors}\label{sec:formal}

Let us first explain, at a formal level, how the correctors are constructed, why the classical correctors cannot be used in a bounded domain, and how we modify them. Recall that $u=u^\eps_N$ is the solution of the Galerkin equation \eqref{eq:galerkin} and that the generator of $(u,m^\eps)$ is given by \eqref{eq:generator}.

\emph{Classical correctors.} Let $\Phi$ be a test function of $u$. In $\Lc^\eps\Phi=-\eps^{-1}D_u\Phi[\sff]-D_u\Phi[\sfs]$, the first term is singular. Following~\cite{PSV}, we look for a perturbed test function $\Phi+\eps\Phi_1$, where $\Phi_1=\Phi_1(u,n)$ also depends on the fast variable. By \eqref{eq:generator},
\[
 \Lc^\eps(\Phi+\eps\Phi_1)=\textcolor{blue}{\eps^{-1}M\Phi_1-\eps^{-1}D_u\Phi[\sff]}-D_u\Phi[\sfs]-D_u\Phi_1[\sff]-\eps D_u\Phi_1[\sfs],
\]
where the singular terms are written in blue. We choose $\Phi_1$ such that they cancel, that is
\[
 M\Phi_1=D_u\Phi[\sff],
\]
where $M$ acts on the variable $n$, $u$ being a parameter. Since $\sff$ is affine in $u$ and $\Phi$ is either $(u,h)$, $(u,h)^2$, $|u|^2$ or $|u|^4$, $D_u\Phi[\sff]$ is a sum of terms of the form $\Gamma(u)c(n)$ and $\Gamma(u)(u,\psi(n))$, where $\Gamma$ is a polynomial, $c$ is a real-valued and $\psi$ an $H$-valued centered observable of the fast variable. For instance, if $\Phi(u)=(u,h)$ with $h$ a test vector, then by \ref{F4}
\[
 D_u\Phi[\sff]=(\sff,h)=(\nN,Ah)-((\nN\cdot\nabla)h,u_s)+(u,F_h(n)),\qquad F_h(n)=-B(\nN,h),
\]
and if $\Phi(u)=|u|^2$, then $D_u\Phi[\sff]=2(u,A\nN)+2(u,B(\nN,u_s))$. The cancellation in \ref{F4} is crucial for this latter test function. The classical corrector associated with $\Gamma(u)(u,\psi(n))$ is
\[
 \Gamma(u)(u,R_0\psi(n)),\qquad R_0\psi:=M^{-1}\psi=-\int_0^\infty\mathsf P_t\psi\dd t .
\]
With this choice of $\Phi_1$,
\[
 \Lc^\eps(\Phi+\eps\Phi_1)=-D_u\Phi[\sfs]-D_u\Phi_1[\sff]-\eps D_u\Phi_1[\sfs].
\]
The term $D_u\Phi_1[\sff]$ is of order one, but it depends on the fast variable $n$. It is replaced by its average with respect to $\eta$ with the help of a second corrector $\eps^2\Phi_2$. By \eqref{eq:generator},
\[
 \Lc^\eps(\Phi+\eps\Phi_1+\eps^2\Phi_2)=-D_u\Phi[\sfs]\textcolor{blue}{-D_u\Phi_1[\sff]+M\Phi_2}-\eps D_u\Phi_1[\sfs]-\eps D_u\Phi_2[\sff]-\eps^2D_u\Phi_2[\sfs],
\]
and we choose $\Phi_2$ such that
\[
 M\Phi_2=D_u\Phi_1[\sff]-\E_\eta D_u\Phi_1[\sff],
\]
which is possible since the right-hand side is centered. Then
\[
 \Lc^\eps(\Phi+\eps\Phi_1+\eps^2\Phi_2)=\overline{\Lc}\Phi-\eps D_u\Phi_1[\sfs]-\eps D_u\Phi_2[\sff]-\eps^2D_u\Phi_2[\sfs],
\]
where
\[
 \overline{\Lc}\Phi(u):=-D_u\Phi(u)[\sfs]-\E_\eta D_u\Phi_1(u,\cdot)[\sff].
\]
The last three terms are formally of order $\eps$, and the averaged operator $\overline{\Lc}$ is the generator of the limit equation.

\emph{Why the classical correctors cannot be used.} The perturbed generator contains $-\eps D_u\Phi_1[\sfs]$, and in particular the term
\[
 -\eps\Gamma(u)(Au,R_0\psi(n))
\]
coming from the Stokes part of the slow drift $\sfs=Au+\PN B(u,v)$. On the torus, or if $R_0\psi$ takes values in $D(A^{1/2})$, it equals $-\eps\Gamma(u)(A^{1/2}u,A^{1/2}R_0\psi)$ and it is bounded by $C\eps\|u\|\,\|R_0\psi\|$. In a bounded domain, this is not the case. Indeed, if for example $\psi(n)=F_h(n)=-P((\nN\cdot\nabla)h)$, then $\psi(n)$ is divergence free and belongs to $H^1$, but the Leray projector does not preserve the boundary condition, so that in general $\psi(n)\notin V$. By \ref{F1}, $\psi(n)$, and thus $R_0\psi(n)$, only belongs to $D(A^\theta)$ for $\theta<1/4$. The term $(Au,R_0\psi)$ then cannot be controlled by the energy since it requires a control on $u$ in $D(A^s)$ for $s>3/4$, whereas the energy gives a control for $s=1/2$. At the level of the Galerkin approximations, $|(Au,\PN R_0\psi)|\le\|u\|\,\|\PN R_0\psi\|$, and the right-hand side is not bounded uniformly in $N$. The same problem occurs for $\Phi=|u|^2$, with $\psi(n)=2B(\nN,u_s)$.

\emph{Modified correctors.} The bad term comes from the action of the Stokes operator on the factor $u$ which is paired with the coefficient of the corrector. Since $A$ is self-adjoint, $(Au,g)=(u,Ag)$, so that this term may be seen as the action of $A$ on the coefficient $g$. Therefore, for a corrector $\Gamma(u)(u,g(n))$, the generator \eqref{eq:generator} contains the two terms
\[
 \eps^{-2}\Gamma(u)(u,Mg)-\Gamma(u)(Au,g)=\eps^{-2}\Gamma(u)(u,(M-\eps^2A)g).
\]
This suggests to keep a corrector of the same form, linear in the factor $u$, but to choose its coefficient as the solution of
\[
 (M-\eps^2A)g=\PN\psi,\qquad\text{that is}\qquad g=\PN R_\eps\psi,\quad R_\eps:=(M-\eps^2A)^{-1}.
\]
The generator of the fast process and the Stokes part of the slow drift are thus inverted together. With this choice, $\eps^{-2}M$ applied to the corrector $\eps\Gamma(u)(u,g)$ still produces $\eps^{-1}\Gamma(u)(u,\psi)$, which compensates the singular term, and the Stokes term $-\eps\Gamma(u)(Au,g)$ disappears exactly. The operator $A$ does not act on $u$: it acts on the values of $\psi$ and $g$, and the function of $u$ is not changed. Only the Stokes part of the slow drift can be treated in this way, since it is linear and self-adjoint. The nonlinear part $\PN B(u,v)$ and the action of $A$ on the factor $\Gamma(u)$ remain; they give the terms $(u,Ah)$ or $\|u\|^2$, which are harmless. The price to pay is that $R_\eps\psi$ satisfies the boundary condition but develops a boundary layer of thickness $\eps$, the distance over which viscosity acts during one correlation time: its $V$-norm is of order $\eps^{-1+2\theta}$ (Lemma~\ref{lem:R}).

Quadratic correctors $\langle u\otimes u,\mathcal K(n)\rangle$ are treated in the same way, since $(Au,\mathcal K(n)u)+(u,\mathcal K(n)Au)=\langle u\otimes u,A\mathcal K(n)+\mathcal K(n)A\rangle$: one inverts $M-\eps^2\A$ with $\A=A\otimes I+I\otimes A$. For correctors which do not contain $u$, such as $\Gamma(u)c(n)$, or whose coefficient takes values in $E$, no such term appears, and we use $R_0=M^{-1}$.

In the rest of this section, we define the operators $R_\eps$ on a suitable class of observables, prove their main properties, and justify the above computation (Corollary~\ref{cor:cancel}).

\subsection{Temporal decay of observables}
Let $X$ be a separable Hilbert space. A measurable map $\psi:\mathscr E\to X$ with $|\psi(n)|_X\le C\mathrm w_{p'}(n)$ for some $p'\le q$ and $\int \psi\dd\eta=0$ is said to have \emph{temporal decay} $(\varpi,p)$ if
\[
 |\mathsf P_t\psi(n)|_X\le\varpi(t)\,\mathrm w_p(n),\qquad t\ge0,\ n\in\mathscr E,
\]
where $\varpi\ge0$ is bounded and $[\varpi]:=\sup\varpi+\int_0^\infty(1+t)\varpi(t)\dd t<\infty$.

\begin{lemma}[Linear and two-point observables]\label{lem:profiles}
\begin{enumerate}[label=(\roman*)]
\item If $\Lambda:E\to X$ is bounded and linear, then its restriction $\psi=\Lambda$ to $\mathscr E$ has temporal decay $(\|\Lambda\|\varrho_0,r_1)$.
\item Let $(\mathcal B_s)_{s\ge0}$ be bounded bilinear maps $E\times E\to X$, measurable in $s$, such that $\beta:=\int_0^\infty\|\mathcal B_s\|\varrho_2(s)\dd s$ is finite. Then
\[
 \psi(n):=\int_0^\infty \mathcal B_s(n\otimes\E\,\widetilde m(s,n)-\E_\eta[m(0)\otimes m(s)])\dd s
\]
is well defined, centered, and has temporal decay $(\beta\varrho_1,r_2)$.
\end{enumerate}
\end{lemma}

\begin{proof}
(i) Since $\Lambda$ is linear and bounded, $\mathsf P_t\Lambda(n)=\Lambda(\E\,\widetilde m(t,n))$, and the result follows from \eqref{eq:mix1}. Moreover $\int\Lambda\dd\eta=\Lambda(\E_\eta n)=0$.

(ii) By \eqref{eq:mix2} with $t=0$, the integrand is bounded by $\|\mathcal B_s\|\varrho_1(0)\varrho_2(s)\mathrm w_{r_2}(n)$, so that $\psi$ is well defined and $|\psi(n)|_X\le\beta\varrho_1(0)\mathrm w_{r_2}(n)$.

By the Markov property,
\[
 \mathsf P_t\big[n\mapsto n\otimes\E\,\widetilde m(s,n)\big](n)=\E[\widetilde m(t,n)\otimes\widetilde m(t+s,n)],
\]
so that
\[
 \mathsf P_t\psi(n)=\int_0^\infty \mathcal B_s\big(\E[\widetilde m(t,n)\otimes\widetilde m(t+s,n)]-\E_\eta[m(0)\otimes m(s)]\big)\dd s .
\]
Using again \eqref{eq:mix2}, we obtain $|\mathsf P_t\psi(n)|_X\le\beta\varrho_1(t)\mathrm w_{r_2}(n)$.

Finally, by stationarity and the Markov property, $\E_\eta[n\otimes\E\,\widetilde m(s,n)]=\E_\eta[m(0)\otimes m(s)]$, hence $\int \psi\dd\eta=0$.
\end{proof}

All the observables used below are of one of these two types. For instance, if $\Lambda_1,\Lambda_2$ are bounded linear maps on $E$ and $(\Gamma_s)_{s\ge0}$ are bounded linear maps on the range of $\Lambda_2$, then $n\mapsto\Lambda_1(n)\otimes\Lambda_2(L_1(n))$ and $n\mapsto\Lambda_1(n)\otimes\int_0^\infty\Gamma_s\Lambda_2(\E\,\widetilde m(s,n))\dd s$ are, after centering, of type (ii).

\subsection{Definition and properties}
Let $X$ be one of the following spaces, endowed with a non-negative self-adjoint operator $\A$:
\begin{itemize}
\item $X=\R$ or a Hilbert space on which $\A=0$;
\item $X=H$, with $\A=A$;
\item $X=H\otimes H$, the Hilbert--Schmidt operators on $H$, with $\A=A\otimes I+I\otimes A$, that is $\A K=AK+KA$.
\end{itemize}
For $\psi:\mathscr E\to X$, with $X$ one of the spaces above, with temporal decay $(\varpi,p)$ and $\eps\ge0$, set
\begin{equation}\label{eq:Reps}
 R_\eps \psi(n):=-\int_0^\infty e^{-t\eps^2\A}\mathsf P_t\psi(n)\dd t .
\end{equation}
Formally $R_\eps=(M-\eps^2\A)^{-1}$, and $R_0=M^{-1}$. On $H$, $R_\eps$ combines the fast semigroup with the Stokes semigroup at time scale $\eps^2$. For instance, for the Ornstein--Uhlenbeck fields of Example~\ref{ex:OU} and $\psi$ linear, $R_\eps \psi=-(I+\eps^2A)^{-1}\psi$.

\begin{lemma}\label{lem:R}
There exists $C$ depending only on $\theta$ and $S$ such that, for $\psi:\mathscr E\to X$ with temporal decay $(\varpi,p)$ and every $n\in\mathscr E$:
\begin{enumerate}[label=(\roman*)]
\item we have
\[
 |R_\eps \psi(n)|_X\le[\varpi]\,\mathrm w_p(n);
\]
\item if $X=H$ and $\beta\in[0,1)$, then
\[
 |A^\beta R_\eps \psi(n)|\le C_\beta\eps^{-2\beta}[\varpi]\,\mathrm w_p(n);
\]
if $X=H\otimes H$, then
\[
 |(A^{1/2}\otimes I)R_\eps \psi(n)|_{\HS}+|(I\otimes A^{1/2})R_\eps \psi(n)|_{\HS}\le C\eps^{-1}[\varpi]\,\mathrm w_p(n);
\]
\item if $X=H$, $\psi$ takes values in $D(A^\theta)$ and $A^\theta \psi$ has temporal decay $(\varpi_\theta,p)$, then
\[
 \|R_\eps \psi(n)\|\le C\eps^{-1+2\theta}[\varpi_\theta]\,\mathrm w_p(n)
\]
and
\[
 |R_\eps \psi(n)-R_0\psi(n)|\le C\eps^{2\theta}[\varpi_\theta]\,\mathrm w_p(n);
\]
\item for every $N$, let $g=\PN R_\eps \psi$ (respectively $(\PN\otimes\PN)R_\eps \psi$, respectively $R_\eps \psi$ if $\A=0$). Then, for the stationary process $m$,
\[
 g(m(t))-g(m(0))-\int_0^t(\eps^2\A g+\PN \psi)(m(r))\dd r
\]
is an $(\mathcal F_t)$-martingale. In other words $(M-\eps^2\A)g=\PN \psi$.
\end{enumerate}
\end{lemma}

\begin{proof}
(i) is immediate, since $e^{-t\eps^2\A}$ is a contraction.

(ii) By \ref{F2},
\[
 |A^\beta R_\eps \psi(n)|\le C_\beta\int_0^\infty (t\eps^2)^{-\beta}\varpi(t)\dd t\,\mathrm w_p(n),
\]
and
\[
 \int_0^\infty t^{-\beta}\varpi(t)\dd t\le\frac{\sup\varpi}{1-\beta}+\int_1^\infty\varpi(t)\dd t .
\]
The case $H\otimes H$ is identical, since $e^{-\tau\A}=e^{-\tau A}\otimes e^{-\tau A}$.

(iii) We argue similarly, using $A^\theta\mathsf P_t\psi=\mathsf P_tA^\theta\psi$ and
\[
 |A^{1/2}e^{-\tau A}y|\le C\tau^{-1/2+\theta}|A^\theta y|\qquad\text{and}\qquad|(e^{-\tau A}-I)y|\le C\tau^\theta|A^\theta y| .
\]

(iv) The range of $g$ is finite dimensional, thus $\A$ is bounded there and commutes with $\mathsf P_r$. By the semigroup property and Fubini's theorem, justified by the temporal decay,
\[
 \mathsf P_rg=e^{r\eps^2\A}(g+\int_0^re^{-s\eps^2\A}\mathsf P_s\PN \psi\dd s)\quad\text{pointwise},
\]
which is the variation-of-constants formula for $\phi'=\eps^2\A\phi+\mathsf P_r\PN \psi$, $\phi(0)=g$. Hence
\[
 \mathsf P_rg-g=\int_0^r\mathsf P_s(\eps^2\A g+\PN \psi)\dd s .
\]
The martingale property follows from the Markov property, since $|g|+|\eps^2\A g|+|\PN \psi|\le C_{N,\eps}\mathrm w_{p\vee p'}$ and $\E \mathrm w_{p\vee p'}(m(t))<\infty$.
\end{proof}

We can now justify the formal computation of Section~\ref{sec:formal}.

\begin{corollary}\label{cor:cancel}
Let $N\ge1$, let $u=u^\eps_N$ be the solution of the Galerkin equation \eqref{eq:galerkin}, and let $\Gamma:H_N\to\R$ be a polynomial.
\begin{enumerate}[label=(\roman*)]
\item Let $\psi:\mathscr E\to H$ have temporal decay $(\varpi,p)$, let $g=\PN R_\eps \psi$ and $\Phi(u,n)=\Gamma(u)(u,g(n))$. Then $\Phi(u_t,m^\eps_t)-\Phi(u_0,m^\eps_0)-\int_0^t\Lc^\eps_\sharp\Phi\dd r$ is a local martingale, where
\[
 \Lc^\eps_\sharp\Phi=\eps^{-2}\Gamma(u)(u,\psi)-\eps^{-1}D_u\Phi[\sff]-D\Gamma(u)[\sfs]\,(u,g)-\Gamma(u)(\PN B(u,v),g).
\]
\item Let $\psi:\mathscr E\to H\otimes H$ have temporal decay $(\varpi,p)$, let $\mathcal K=(\PN\otimes\PN)R_\eps \psi$ and $\Phi(u,n)=\Gamma(u)\langle u\otimes u,\mathcal K(n)\rangle$. Then the same holds with
\[
 \Lc^\eps_\sharp\Phi=\eps^{-2}\Gamma(u)\langle u\otimes u,\psi\rangle-\eps^{-1}D_u\Phi[\sff]-D\Gamma(u)[\sfs]\,\langle u\otimes u,\mathcal K\rangle-\Gamma(u)(\PN B(u,v),(\mathcal K+\mathcal K^*)u).
\]
\end{enumerate}
\end{corollary}

\begin{proof}
Recall that $t\mapsto u_t$ is a continuous adapted process with values in the finite-dimensional space $H_N$; it is absolutely continuous, with derivative $\partial_tu=-\eps^{-1}\sff(u,m^\eps)-\sfs(u)$, which is c\`adl\`ag. Moreover $g(m^\eps_t)$ is a semimartingale by Lemma~\ref{lem:R}(iv) and the time change $t\mapsto t/\eps^2$:
\[
 \dd g(m^\eps_t)=\eps^{-2}(\eps^2Ag+\PN \psi)(m^\eps_t)\dd t+\dd\mathcal M_t,
\]
where $\mathcal M$ is a martingale for the filtration $(\mathcal F_{t/\eps^2})_{t\ge0}$.
The product rule has no bracket term, because $u$ is continuous with finite variation. Therefore, in $\dd\Phi$, the term $-\Gamma(u)(Au,g)\dd t$ coming from $\dd u$, which is equal to $-\Gamma(u)(u,Ag)\dd t$, cancels with the term $\Gamma(u)(u,\eps^{-2}\eps^2Ag)\dd t$ coming from $\dd g$. The martingale term $\int\Gamma(u)(u,\dd\mathcal M)$ is a local martingale. In case (ii), the terms $-\Gamma(u)((Au,\mathcal Ku)+(u,\mathcal KAu))$ coming from $\dd u$ cancel in the same way with $\Gamma(u)\langle u\otimes u,A\mathcal K+\mathcal KA\rangle$, which comes from $\dd\mathcal K(m^\eps_t)$.
\end{proof}

As seen in Section~\ref{sec:formal}, we have overcome the problem of the classical correctors. Indeed, compared with the formal generator \eqref{eq:generator}, the terms
\[
 -\Gamma(u)(Au,g)\quad\text{in case (i)}\qquad\text{and}\qquad-\Gamma(u)(Au,(\mathcal K+\mathcal K^*)u)\quad\text{in case (ii)}
\]
have disappeared. The Stokes operator is thus moved from $u$, which is only in $V$, onto the corrector, which lies in $D(A)$. What we lose is the bound $\|R_\eps \psi\|\lesssim\eps^{-1+2\theta}$ of Lemma~\ref{lem:R}(iii), which reflects a boundary layer of thickness $\eps$. On the torus, or for coefficients $\psi$ with values in $D(A^{1/2})$, one may take $\eps=0$ in $R_\eps$. For correctors with $\A=0$ (scalar or $E$-valued), no cancellation occurs and \eqref{eq:generator} is used as it stands.

\section{Uniform energy estimates}\label{sec:energy}

We now derive bounds which are uniform in $\eps$. To this end, we apply the perturbed test function method with modified correctors to the energy.

In dimension three, the energy equality is not known to hold for solutions of \eqref{eq:eps}, so the martingale identities below cannot be applied to $|u|^2$ or $|u|^4$ along $u^\eps$. We therefore perform all the computations of this section and of Sections~\ref{sec:ptf}--\ref{sec:tight} on $u_N^\eps$, with constants independent of $N$ and $\eps\in(0,1)$. This gives directly Theorem~\ref{thm:main}(ii), for any $N(\eps)\to\infty$, and part (iii) follows by a diagonal argument, which we now explain.

For fixed $\eps$, the solution $u^\eps$ of \eqref{eq:eps} has been obtained in Section~\ref{sec:existence} as the limit in law of a subsequence $(u^\eps_{N_j})_{j\ge1}$ of the Galerkin approximations, which may depend on $\eps$; in dimension three, we do not know whether the whole sequence converges. The convergence of Section~\ref{sec:existence} holds jointly with $m^\eps$, that is $\Law(u^\eps_{N_j},m^\eps)\to\Law(u^\eps,m^\eps)$ in $\mathcal X_T\times D([0,T];E)$. Moreover, the map $x\mapsto\eps^{-1}\int_0^\cdot x(s)\dd s$ is continuous from $D([0,T];E)$, endowed with the Skorokhod topology, into $C([0,T];E)$. Indeed, if $x_n\to x$ in $D([0,T];E)$, then $\sup_n\sup_{s\le T}\|x_n(s)\|_E<\infty$ and $x_n(s)\to x(s)$ for all $s$ outside the countable set of discontinuity points of $x$~\cite[Section~12]{Billingsley1999}, so that $\int_0^T\|x_n(s)-x(s)\|_E\dd s\to0$ by dominated convergence. Since $z^\eps$ is the image of $m^\eps$ by this map, we obtain $\Law(u^\eps_{N_j},z^\eps)\to\Law(u^\eps,z^\eps)$ as $j\to\infty$. Let $d_{\mathrm{LP}}$ be the L\'evy--Prokhorov metric, which metrizes the weak convergence of probability measures on $\mathcal X_T\times\mathcal Z_T$ (with $\varsigma\ge\gamma$). For each $\eps$, we may choose an index $N(\eps)$ in this subsequence such that
\[
 N(\eps)\ge\eps^{-1}\qquad\text{and}\qquad d_{\mathrm{LP}}(\Law(u^\eps_{N(\eps)},z^\eps),\Law(u^\eps,z^\eps))\le\eps .
\]
Then, for every sequence $\eps_k\to0$, the sequences $\Law(u^{\eps_k},z^{\eps_k})$ and $\Law(u^{\eps_k}_{N(\eps_k)},z^{\eps_k})$ have the same limit points, and (iii) follows from (ii).

In dimension two, by Proposition~\ref{prop:2D}, the whole Galerkin sequence converges almost surely to $u^\eps$, and no justification is necessary.
Let us now describe the strategy, which follows the formal construction of Section~\ref{sec:formal}. For $\Phi(u)=|u|^2$, a first corrector $\eps\Phi_1$ removes the singular terms of order $\eps^{-1}$; it contains the modified corrector $\Psi_E$. The terms of order one which depend on the fast variable are quadratic in $n$, hence not bounded along $m^\eps$, and a second corrector $\eps^2\Phi_2$ replaces them by their average. The remainders are of order $\eps^{2\theta}\mathrm w_p(n)$: the loss of $\eps^{1-2\theta}$ comes from the boundary layer of the modified correctors (Lemma~\ref{lem:R}(iii)), and the weight $\mathrm w_p(m^\eps(t))$ is bounded by $2\eps^{-\alpha p}$ up to the stopping time $\tau^\eps$. We then repeat the construction for $|u|^4$ (Lemma~\ref{lem:quartic}), which gives the quartic bound used in Section~\ref{sec:ptf} to control the remainders.

Throughout this section, $\Phi(u)=|u|^2$, and all constants are independent of $N$ and $\eps\in(0,1)$. We use repeatedly the elementary bounds, valid for $u,w\in H_N$ and $n\in E$, which follow from \ref{F3}--\ref{F4}:
\begin{gather}
 |(\sff,w)|\le C\|n\|_E(1+\|u\|)|w|,\qquad|(\sff,h)|\le C\|n\|_E(1+|u|),\qquad|(\sfs,h)|\le C(1+|u|^2),\label{eq:elem}\\
 |(\PN B(u,v),w)|=|((u\cdot\nabla)w,v)|\le C(|u|^{1/4}\|u\|^{3/4}+|u|^{1/2}\|u\|^{3/2})\|w\| .\nonumber
\end{gather}

\paragraph{First corrector.} Let $F_E(n):=B(\nN,u_s)$. It is linear in $n$, with values in $H\cap H^1\subset D(A^\theta)$, so $F_E$ and $A^\theta F_E$ have temporal decay $(C\varrho_0,r_1)$ by Lemma~\ref{lem:profiles}(i). Set $\Psi_E:=\PN R_\eps F_E$; then $|\Psi_E|\le C\mathrm w_{r_1}$ and $\|\Psi_E\|\le C\eps^{-1+2\theta}\mathrm w_{r_1}$ by Lemma~\ref{lem:R}. Let
\[
 \Phi_1(u,n)=2(u,A\PN L_1(n))+2(u,\Psi_E(n)),\qquad|\Phi_1|\le C|u|\mathrm w_{r_1}(n).
\]
Recall that $(B(\nN,u),u)=(B(u,u),u)=0$ by \ref{F4}, so that $\Lc^\eps\Phi=-2\eps^{-1}(u,A\nN+B(\nN,u_s))-2\|u\|^2-2(B(u,u_s),u)$. Since $M(u,A\PN L_1)=(u,A\nN)$, Corollary~\ref{cor:cancel} gives
\[
 \Lc^\eps_\sharp(\Phi+\eps\Phi_1)=-2\|u\|^2-2(B(u,u_s),u)-\Theta(u,n)-\eps S_1(u,n),
\]
where the terms of order $\eps^{-1}$ have canceled, and
\[
 \Theta:=D_u\Phi_1[\sff]=2(\sff,A\PN L_1+\Psi_E),\qquad S_1:=2(u,A^2\PN L_1)+2(\PN B(u,v),A\PN L_1+\Psi_E).
\]

\paragraph{Second corrector.} The order-one term $\Theta$ is quadratic in the fast variable, hence not bounded along $m^\eps$. As in~\cite{BarrueDebusscheTusseau2024}, it is replaced by its average with a second corrector. By \ref{F4}, $\Theta(u,n)=q_E(n)+(u,g_E(n))$ with
\begin{gather*}
 q_E(n)=2(A\nN+B(\nN,u_s),A\PN L_1(n)+\Psi_E(n)),\\
 g_E(n)=-2\PN P((\nN\cdot\nabla)(A\PN L_1(n)+\Psi_E(n))).
\end{gather*}
Since $L_1(n)=-\int_0^\infty\E\,\widetilde m(s,n)\dd s$ and $\Psi_E(n)=-\int_0^\infty e^{-s\eps^2A}\PN B(\PN\E\,\widetilde m(s,n),u_s)\dd s$, the centered observables $q_E-\E_\eta q_E$ and $g_E-\E_\eta g_E$ are of type (ii) in Lemma~\ref{lem:profiles}. For $q_E$, the bilinear maps are bounded uniformly in $s$. For $g_E$, they are bounded by $C(1+(s\eps^2)^{-1/2+\theta})$, because $|(n\cdot\nabla)e^{-s\eps^2A}y|\le C\|n\|_E\|e^{-s\eps^2A}y\|$ and \ref{F1}--\ref{F2} apply. Hence, with $\chi_E:=R_0(q_E-\E_\eta q_E)$ and $\Xi_E:=\PN R_\eps(g_E-\E_\eta g_E)$,
\[
 |\chi_E|\le C\mathrm w_{r_2},\qquad|\Xi_E|\le C\eps^{-1+2\theta}\mathrm w_{r_2},\qquad\|\Xi_E\|\le C\eps^{-2+2\theta}\mathrm w_{r_2},
\]
and we set $\Phi_2:=\chi_E+(u,\Xi_E)$. The average $\overline\Theta(u):=\E_\eta\Theta(u,\cdot)$ is deterministic. Moreover $(u,\E_\eta g_E)=2\E_\eta((\nN\cdot\nabla)u,A\PN L_1+\Psi_E)$, so that
\begin{equation}\label{eq:Qbar}
 |\overline\Theta(u)|\le C\E_\eta[\|n\|_E\mathrm w_{r_1}(n)](1+\|u\|)\le\tfrac12\|u\|^2+C .
\end{equation}
\paragraph{Remainders.} By Corollary~\ref{cor:cancel}, with $\Phi_E^\eps:=\Phi+\eps\Phi_1+\eps^2\Phi_2$,
\begin{equation}\label{eq:energy-identity}
 \Lc^\eps_\sharp\Phi_E^\eps=-2\|u\|^2-2(B(u,u_s),u)-\overline\Theta(u)-\mathcal R^E_\eps,\qquad
 \mathcal R^E_\eps=\eps S_1+\eps(\sff,\Xi_E)+\eps^2(\PN B(u,v),\Xi_E).
\end{equation}
By \eqref{eq:elem} and the bounds above,
\[
 |\mathcal R^E_\eps|\le C\eps^{2\theta}\mathrm w_p(n)[(1+|u|)^2+(1+\|u\|)+|u|^{1/4}\|u\|^{3/4}+|u|^{1/2}\|u\|^{3/2}].
\]
For $t<\tau^\eps$, Lemma~\ref{lem:stop} gives $\mathrm w_p(m^\eps(t))\le2\eps^{-\alpha p}$ and $\eps^{2\theta-\alpha p}\le1$ by \eqref{eq:theta-alpha}. Young's inequality then yields, for $t<\tau^\eps$,
\begin{equation}\label{eq:energy-gen}
 \Lc^\eps_\sharp\Phi_E^\eps(u,m^\eps(t))\le-\|u\|^2+C(1+|u|^2).
\end{equation}
The boundary correctors satisfy $\eps|\Phi_1|+\eps^2|\Phi_2|\le\frac12|u|^2+C\eps^2\mathrm w_{2p}(n)$.

\begin{remark}[Boundary terms at the stopping time]\label{rem:endpoint}
The corrector terms evaluated at the final time $t\wedge\tau^\eps$ involve $m^\eps(t\wedge\tau^\eps)$. As in~\cite[Section~6]{BarrueDebusscheTusseau2024}, if $\tau^\eps>0$ and $m$ has continuous paths (as in Example~\ref{ex:OU}), then $\|m^\eps(t\wedge\tau^\eps)\|_E\le\eps^{-\alpha}$, and these terms are bounded by $C\eps^{1-\alpha p}(1+|u|^2)$ pathwise. In general, $m$ is only c\`adl\`ag and may jump above $\eps^{-\alpha}$ exactly at $\tau^\eps$; when $\tau^\eps=0$, the value $m^\eps(0)$ is arbitrary. In both situations we use \eqref{eq:supmoment} instead, which gives $\eps^{c}\E\sup_{t\le T}\mathrm w_{p'}(m^\eps(t))\le C\eps^{c-2p'/q}$ for every weight occurring below. This is how all boundary terms are handled in what follows.
\end{remark}

\begin{lemma}[Correctors for the quartic energy]\label{lem:quartic}
Let $\Psi_1:=2|u|^2\Phi_1$ and $\Psi_2:=2|u|^2\Phi_2+8\langle u\otimes u,\mathcal K_E\rangle$, where $\mathcal K_E:=(\PN\otimes\PN)R_\eps(k_E-\E_\eta k_E)$ and $k_E:=(A\nN+B(\nN,u_s))\otimes(A\PN L_1+\Psi_E)$. Then
\[
 \Lc^\eps_\sharp(|u|^4+\eps\Psi_1+\eps^2\Psi_2)=-4|u|^2(u,\sfs)-\overline\Theta_4(u)-\mathcal R^{(4)}_\eps,
\]
where $\overline\Theta_4:=\E_\eta\Theta_4$, $\Theta_4:=2|u|^2\Theta+8(u,\sff)(u,A\PN L_1+\Psi_E)$, and, for $t<\tau^\eps$,
\[
 |\mathcal R^{(4)}_\eps(u,m^\eps(t))|\le\tfrac12|u|^2\|u\|^2+C(1+|u|^4+\|u\|^2).
\]
\end{lemma}

\begin{proof}
Note that $(u,\sff)=(u,A\nN+B(\nN,u_s))$ and $(u,\sfs)=\|u\|^2+(B(u,u_s),u)$, by \ref{F4}, and that $D(|u|^4)[k]=4|u|^2(u,k)$ and
\[
 \Lc^\eps|u|^4=-4\eps^{-1}|u|^2(u,\sff)-4|u|^2(u,\sfs) .
\]

\emph{First corrector.} By Corollary~\ref{cor:cancel}, $\eps^{-2}M$ applied to $\eps\Psi_1$ gives $4\eps^{-1}|u|^2(u,\sff)$, which cancels the singular term. The remaining order-one term is
\[
 -D_u\Psi_1[\sff]=-\Theta_4,\qquad\Theta_4:=2|u|^2\Theta(u,n)+4(u,\sff)\,\Phi_1=2|u|^2\Theta+8(u,\sff)\,(u,A\PN L_1+\Psi_E).
\]

\emph{Second corrector.} By the structure of $\Theta$, the fluctuation of $\Theta_4$ is $2|u|^2[(q_E-\E_\eta q_E)+(u,g_E-\E_\eta g_E)]+8\langle u\otimes u,k_E-\E_\eta k_E\rangle$, where $k_E$ is a two-point observable in $H\otimes H$. This explains the choice of $\Psi_2$, and $|\mathcal K_E|_{\HS}\le C\mathrm w_{r_2}$ and $\|\mathcal K_E\|_{V,\HS}\le C\eps^{-1}\mathrm w_{r_2}$ by Lemma~\ref{lem:R}. By Corollary~\ref{cor:cancel},
\[
 \Lc^\eps_\sharp(|u|^4+\eps\Psi_1+\eps^2\Psi_2)=-4|u|^2(u,\sfs)-\overline\Theta_4(u)-\mathcal R^{(4)}_\eps,
\]
with $\overline\Theta_4:=\E_\eta\Theta_4$ and $\mathcal R^{(4)}_\eps=U_1+\dots+U_8$, where
\begin{align*}
 U_1&=4\eps(u,\sfs)\Phi_1, & U_2&=2\eps|u|^2S_1, & U_3&=4\eps(u,\sff)\Phi_2,\\
 U_4&=2\eps|u|^2(\sff,\Xi_E), & U_5&=8\eps(\sff,(\mathcal K_E+\mathcal K_E^*)u), & U_6&=4\eps^2(u,\sfs)\Phi_2,\\
 U_7&=2\eps^2|u|^2(\PN B(u,v),\Xi_E), & U_8&=8\eps^2(\PN B(u,v),(\mathcal K_E+\mathcal K_E^*)u).
\end{align*}
Let $t<\tau^\eps$ and $\mathsf v=\mathrm w_p(m^\eps(t))$, so that $\eps^{2\theta}\mathsf v\le2$ and $\eps \mathsf v\le2$. Using $|(u,\sfs)|\le\|u\|^2+C|u|^2$, $|(u,\sff)|\le C\|n\|_E|u|$, $|\Phi_1|\le C|u|\mathrm w_{r_1}$, $|\Phi_2|\le C\mathrm w_{r_2}(1+\eps^{-1+2\theta}|u|)$, the bounds on $\Xi_E$, $\mathcal K_E$ and \eqref{eq:elem}, we get:
\begin{align*}
 |U_1|+|U_6|&\le C(\|u\|^2+|u|^2)(|u|+1)\le\tfrac18|u|^2\|u\|^2+C(1+|u|^4+\|u\|^2),\\
 |U_2|+|U_7|&\le C|u|^2(1+|u|^2+|u|^{1/4}\|u\|^{3/4}+|u|^{1/2}\|u\|^{3/2})\le\tfrac18|u|^2\|u\|^2+C(1+|u|^4),\\
 |U_3|&\le C(|u|+|u|^2),\\
 |U_4|&\le C|u|^2(1+\|u\|)\le\tfrac18|u|^2\|u\|^2+C(1+|u|^2),\\
 |U_5|&\le C|u|(1+\|u\|)\le C(1+|u|^2+\|u\|^2),\\
 |U_8|&\le C|u|(|u|^{1/4}\|u\|^{3/4}+|u|^{1/2}\|u\|^{3/2})\le\tfrac18|u|^2\|u\|^2+C(1+|u|^2+\|u\|^2).
\end{align*}
In $U_1$ and $U_6$ we used $\|u\|^2|u|\le\delta|u|^2\|u\|^2+C_\delta\|u\|^2$; in $U_2$, $U_7$ and $U_8$, Young's inequality with exponents adapted to the powers of $\|u\|$, which are smaller than $2$. Hence $|\mathcal R^{(4)}_\eps|\le\frac12|u|^2\|u\|^2+C(1+|u|^4+\|u\|^2)$ for $t<\tau^\eps$.

\end{proof}

\begin{proposition}[Uniform stopped energy bounds]\label{prop:energy}
For every $T>0$ there is $C_T$, independent of $N$ and $\eps\in(0,1)$, such that for all stopping times $\tau\le\tau^\eps$
\begin{gather*}
 \E|u(\tau)|^2+\E\int_0^\tau\|u\|^2\dd r\le C_T,\qquad
 \E|u(\tau)|^4+\E\int_0^\tau|u|^2\|u\|^2\dd r\le C_T,\\
 \E\sup_{t\le\tau^\eps}|u(t)|^2\le3C_T^{1/2}.
\end{gather*}
Moreover $\Pp(\tau^\eps<T)\to0$ as $\eps\to0$.
\end{proposition}

\begin{proof}
\emph{Quadratic bound.} Let $\tau_R=\inf\{t:|u(t)|\ge R\}$. We apply Corollary~\ref{cor:cancel} to $\Phi_E^\eps$ between $0$ and $\tau\wedge t\wedge\tau_R$; the stopped local martingale is a martingale, since $|u|\le R$ and $\E\sup_{t\le T}\mathrm w_{2p}(m^\eps(t))<\infty$. We take expectations and use \eqref{eq:energy-gen} for $r<\tau^\eps$ and the bound on the boundary correctors. With \eqref{eq:supmoment}, $\eps^2\E\sup_t\mathrm w_{2p}(m^\eps(t))\le C\eps^{2-4p/q}\le C$, so that
\[
 \tfrac12\E|u(\tau\wedge t\wedge\tau_R)|^2+\E\int_0^{\tau\wedge t\wedge\tau_R}\|u\|^2\dd r\le|u_0|^2+C+C\int_0^t\E|u(\tau\wedge r\wedge\tau_R)|^2\dd r .
\]
Gronwall's lemma and $R\to\infty$ give the bound.

\emph{Quartic bound.} By \eqref{eq:Qbar} and $|\E_\eta[(u,\sff)\,(u,A\PN L_1+\Psi_E)]|\le C|u|^2\E_\eta[\|n\|_E\mathrm w_{r_1}(n)]$, we have $\overline\Theta_4\ge-|u|^2\|u\|^2-C(1+|u|^4)$. Finally $-4|u|^2(u,\sfs)\le-4|u|^2\|u\|^2+C|u|^4$. Altogether, by Lemma~\ref{lem:quartic}, for $t<\tau^\eps$,
\[
 \Lc^\eps_\sharp(|u|^4+\eps\Psi_1+\eps^2\Psi_2)\le-\tfrac52|u|^2\|u\|^2+C(1+|u|^4+\|u\|^2),
\]
and $\E\int_0^\tau\|u\|^2\dd r\le C_T$ by the quadratic bound.

The boundary correctors satisfy
\[
 \eps|\Psi_1|+\eps^2|\Psi_2|\le\tfrac12|u|^4+C\eps^4\mathrm w_{4p}(n),
\]
and $\eps^4\E\sup_t\mathrm w_{4p}(m^\eps(t))\le C$ by \eqref{eq:supmoment} and $4p<q$. Localizing as in the quadratic case and applying Gronwall's lemma gives the quartic bound.

\emph{Supremum.} The process $X_t=|u(t\wedge\tau^\eps)|^4$ is continuous and adapted, and $\E X_\tau\le C_T$ for every stopping time $\tau\le T$. Lenglart's inequality~\cite{Lenglart1977} gives $\E[(\sup_{t\le T}X_t)^{1/2}]\le3C_T^{1/2}$. The last claim is Lemma~\ref{lem:stop}.
\end{proof}

The remainders of Section~\ref{sec:ptf} will be controlled by the quantity $\Pol(u):=(1+|u|^2)(1+\|u\|^2)$. By Proposition~\ref{prop:energy}, $\sup_{N,\eps}\E\int_0^{\tau^\eps}\Pol(u)\dd t<\infty$.

\section{Perturbed test functions for cylindrical observables}\label{sec:ptf}

Let us now construct the correctors for the observables which are used to prove tightness and to identify the limit. We follow again the three steps of Section~\ref{sec:formal}, now for $\Phi(u)=(u,h)^\ell$, $\ell\in\{1,2\}$. The first corrector contains the modified corrector $\Psi=\PN R_\eps F_h$ and removes the singular terms. The second corrector replaces the order-one terms depending on $n$ by their average, which defines the averaged operator $\Lav_\eps$. The remainders are estimated in Lemma~\ref{lem:remainders}; they are of order $\eps^{2\theta}$ and are controlled by the quantity $\Pol(u)$ of Section~\ref{sec:energy}. In Section~\ref{sec:z}, the same construction is applied to the integrated fast velocity. Fix $h\in H_{N_0}$, $\ell\in\{1,2\}$, and let $\Phi(u)=(u,h)^\ell$. For $u\in H_N$, by \ref{F4},
\begin{equation}\label{eq:fh}
 (\sff,h)=\mathsf d(n)+(u,F_h(n)),\qquad\mathsf d(n):=(\nN,Ah)-((\nN\cdot\nabla)h,u_s),\qquad F_h(n):=-B(\nN,h).
\end{equation}
Recall that the maps $\mathsf d$ and $F_h$ are linear, and that $F_h$ takes values in $H\cap H^1\subset D(A^\theta)$.

\subsection{Construction of the correctors}

\emph{First corrector.} Set
\[
 \mathsf c(n):=R_0\mathsf d(n)=\mathsf d(L_1(n)),\qquad\Psi:=\PN R_\eps F_h,\qquad\Phi_1:=\ell(u,h)^{\ell-1}\big[\mathsf c(n)+(u,\Psi(n))\big].
\]

By Lemmas~\ref{lem:profiles}(i) and~\ref{lem:R}, $|\mathsf c|+|\Psi|\le C\mathrm w_{r_1}$, $\|\Psi\|\le C\eps^{-1+2\theta}\mathrm w_{r_1}$ and $|\Psi-\Psi_0|\le C\eps^{2\theta}\mathrm w_{r_1}$, where $\Psi_0:=\PN R_0F_h=-\PN B(\PN L_1(n),h)$. By Corollary~\ref{cor:cancel}, the terms of order $\eps^{-1}$ in $\Lc^\eps_\sharp(\Phi+\eps\Phi_1)$ are $\eps^{-1}\ell(u,h)^{\ell-1}[\mathsf d+(u,F_h)]-\eps^{-1}D\Phi[\sff]=0$.

\emph{Order one.} The $n$-dependent term of order one is
\begin{align}
 Y(u,n):=D_u\Phi_1[\sff]&=\ell(\ell-1)(u,h)^{\ell-2}(\mathsf d+(u,F_h))(\mathsf c+(u,\Psi))\label{eq:Y}\\
 &\quad+\ell(u,h)^{\ell-1}[(A\nN,\Psi)+((\nN\cdot\nabla)u_s,\Psi)+(u,\mathcal G)],\nonumber
\end{align}
where, by \ref{F4}, $((\nN\cdot\nabla)u,\Psi)=(u,\mathcal G)$ with $\mathcal G(n):=-\PN P((\nN\cdot\nabla)\Psi(n))$. Note that $|(u,\mathcal G)|\le C\|n\|_E\|u\||\Psi|$, so that $|Y|\le C(1+|u|^2)(1+\|u\|)\mathrm w_p(n)$ uniformly in $\eps$. We define the averaged operators
\begin{equation}\label{eq:Leps}
 \Lav_\eps\Phi(u):=-D\Phi(u)[\sfs]-\E_\eta Y(u,\cdot),\qquad\Lav^N_0\Phi(u):=\text{the same with $\Psi_0$ in place of $\Psi$.}
\end{equation}

\emph{Second corrector.} The fluctuation $Z:=Y-\E_\eta Y$ decomposes as follows. For $\ell=1$,
\[
 Z=z_0(n)+(u,g_1(n)),\qquad z_0:=(A\nN+B(\nN,u_s),\Psi)-\E_\eta[\cdots],\qquad g_1:=\mathcal G-\E_\eta \mathcal G .
\]
For $\ell=2$,
\[
 Z=2[(\mathsf d\mathsf c-\E_\eta\mathsf d\mathsf c)+(u,g_2)+\langle u\otimes u,k_2\rangle]+2(u,h)[z_0+(u,g_1)],
\]
with $g_2:=\mathsf d\Psi+\mathsf cF_h-\E_\eta[\cdots]$ and $k_2:=F_h\otimes\Psi-\E_\eta[F_h\otimes\Psi]$. All these are two-point observables of type (ii) in Lemma~\ref{lem:profiles}, since $\mathsf c$ and $\Psi$ are linear functionals of $s\mapsto\E\,\widetilde m(s,n)$. The bilinear maps are bounded uniformly in $s$, except for $g_1$, where they are bounded by $C(1+(s\eps^2)^{-1/2+\theta})$. Set
\begin{gather*}
 \chi_0:=R_0z_0,\qquad\chi_1:=R_0(\mathsf d\mathsf c-\E_\eta\mathsf d\mathsf c),\\
 \Xi_1:=\PN R_\eps g_1,\qquad\Xi_2:=\PN R_\eps g_2,\qquad\mathcal K:=(\PN\otimes\PN)R_\eps k_2,
\end{gather*}
and
\[
 \Phi_2:=\begin{cases}\chi_0+(u,\Xi_1),&\ell=1,\\[0.2em]
 2[\chi_1+(u,\Xi_2)+\langle u\otimes u,\mathcal K\rangle]+2(u,h)[\chi_0+(u,\Xi_1)],&\ell=2.\end{cases}
\]
By Lemma~\ref{lem:R}(i)--(ii),
\begin{gather}
 |\chi_0|+|\chi_1|+|\Xi_2|+|\mathcal K|_{\HS}\le C\mathrm w_{r_2},\qquad\|\Xi_2\|+\|\mathcal K\|_{V,\HS}\le C\eps^{-1}\mathrm w_{r_2},\label{eq:sizes}\\
 |\Xi_1|\le C\eps^{-1+2\theta}\mathrm w_{r_2},\qquad\|\Xi_1\|\le C\eps^{-2+2\theta}\mathrm w_{r_2},\nonumber
\end{gather}
where $\|\mathcal K\|_{V,\HS}:=|(A^{1/2}\otimes I)\mathcal K|_{\HS}+|(I\otimes A^{1/2})\mathcal K|_{\HS}$. Identifying $\mathcal K$ with an operator, $\langle x\otimes y,\mathcal K\rangle=(x,\mathcal Ky)$, we have $D_u\langle u\otimes u,\mathcal K\rangle[k]=(k,(\mathcal K+\mathcal K^*)u)$, $|(\mathcal K+\mathcal K^*)u|\le2|\mathcal K|_{\HS}|u|$ and $\|(\mathcal K+\mathcal K^*)u\|\le\|\mathcal K\|_{V,\HS}|u|$. By Corollary~\ref{cor:cancel} and Lemma~\ref{lem:R}(iv), the order-one part of $\Lc^\eps_\sharp(\eps^2\Phi_2)$ is exactly $Z$.

\subsection{Remainder estimates}

We denote by $\widetilde u(t):=u(t\wedge\tau^\eps)$ the process $u$ stopped at the time $\tau^\eps$ of Lemma~\ref{lem:stop}.

The remainders are estimated in the next lemma. The worst terms are those which involve the $V$-norm of the modified correctors $\Psi$ and $\Xi_1$; through Lemma~\ref{lem:R}(iii), they explain why the remainders are only of order $\eps^{2\theta}$.

\begin{lemma}[Remainders for cylindrical observables]\label{lem:remainders}
Let $\Phi^\eps=\Phi+\eps\Phi_1+\eps^2\Phi_2$. Then $\Lc^\eps_\sharp\Phi^\eps=\Lav_\eps\Phi+\mathcal R_\eps$, with $|\mathcal R_\eps|\le C\eps^{2\theta}\mathrm w_p(n)\Pol(u)$.
\end{lemma}

\begin{proof}
By Corollary~\ref{cor:cancel}, $\mathcal R_\eps=-(T_1+\dots+T_{10})$, where the terms $T_i$ and their bounds are given below. We write them out completely for $\ell=2$; the case $\ell=1$ consists of a subset of the same terms. We set $\mathsf v=\mathrm w_p(n)$, and we use \eqref{eq:elem}, the bounds on $\mathsf c,\Psi$, \eqref{eq:sizes}, and $\|n\|_E\mathrm w_{r_i}(n)\le2\mathrm w_p(n)$. The ten terms are:
\begin{align*}
 T_1&=2\eps(\sfs,h)[\mathsf c+(u,\Psi)], & |T_1|&\le C\eps \mathsf v(1+|u|)^3,\\
 T_2&=2\eps(u,h)(\PN B(u,v),\Psi), & |T_2|&\le C\eps^{2\theta}\mathsf v(|u|^{5/4}\|u\|^{3/4}+|u|^{3/2}\|u\|^{3/2}),\\
 T_3&=2\eps(\sff,\Xi_2), & |T_3|&\le C\eps \mathsf v(1+\|u\|),\\
 T_4&=2\eps(\sff,(\mathcal K+\mathcal K^*)u), & |T_4|&\le C\eps \mathsf v|u|(1+\|u\|),\\
 T_5&=2\eps(\sff,h)[\chi_0+(u,\Xi_1)], & |T_5|&\le C\eps^{2\theta}\mathsf v(1+|u|)^2,\\
 T_6&=2\eps(u,h)(\sff,\Xi_1), & |T_6|&\le C\eps^{2\theta}\mathsf v|u|(1+\|u\|),\\
 T_7&=2\eps^2(\PN B(u,v),\Xi_2), & |T_7|&\le C\eps \mathsf v(|u|^{1/4}\|u\|^{3/4}+|u|^{1/2}\|u\|^{3/2}),\\
 T_8&=2\eps^2(\PN B(u,v),(\mathcal K+\mathcal K^*)u), & |T_8|&\le C\eps \mathsf v(|u|^{5/4}\|u\|^{3/4}+|u|^{3/2}\|u\|^{3/2}),\\
 T_9&=2\eps^2(\sfs,h)[\chi_0+(u,\Xi_1)], & |T_9|&\le C\eps^{1+2\theta}\mathsf v(1+|u|)^3,\\
 T_{10}&=2\eps^2(u,h)(\PN B(u,v),\Xi_1), & |T_{10}|&\le C\eps^{2\theta}\mathsf v(|u|^{5/4}\|u\|^{3/4}+|u|^{3/2}\|u\|^{3/2}),
\end{align*}

Let us detail two of them. For $T_2$, $|(\PN B(u,v),\Psi)|=|((u\cdot\nabla)\Psi,v)|\le|u|_{L^4}|v|_{L^4}\|\Psi\|$, and $\|\Psi\|\le C\eps^{-1+2\theta}\mathrm w_{r_1}$, so the factor $\eps$ becomes $\eps^{2\theta}$; then \ref{F3} gives the stated powers. For $T_{10}$, the same argument applies with $\|\Xi_1\|\le C\eps^{-2+2\theta}\mathrm w_{r_2}$, which absorbs the factor $\eps^2$. This is the most singular term, and it is the reason why the remainders are only of order $\eps^{2\theta}$.

All the other terms involve at most one $V$-norm of a coefficient of size $\eps^{-1}$ (namely $\Xi_2$ or $\mathcal K$) multiplied by $\eps^2$, or an $H$-norm of size $\eps^{-1+2\theta}$ multiplied by $\eps$.

Finally, by Poincar\'e's and Young's inequalities, each of the functions of $u$ on the right is bounded by $C\Pol(u)$: for instance $|u|^3\le C|u|^2\|u\|\le C\Pol(u)$ and $|u|^{3/2}\|u\|^{3/2}\le1+|u|^2\|u\|^2$. Since $\eps\le\eps^{2\theta}$, this proves $|\mathcal R_\eps|\le C\eps^{2\theta}\mathrm w_p(n)\Pol(u)$.
\end{proof}

\begin{proposition}\label{prop:ptf}
Let $\Phi^\eps=\Phi+\eps\Phi_1+\eps^2\Phi_2$. Then:
\begin{enumerate}[label=(\alph*)]
\item $\Phi^\eps(\widetilde u_t,m^\eps_{t\wedge\tau^\eps})-\Phi^\eps(u_0,m^\eps_0)-\int_0^{t\wedge\tau^\eps}\Lc^\eps_\sharp\Phi^\eps\dd r$ is a local martingale; for every $R>0$, it is a martingale when stopped at $\tau_R:=\inf\{t:|\widetilde u(t)|\ge R\}\wedge T$;
\item $|\Phi^\eps-\Phi|\le C\eps \mathrm w_p(n)(1+|u|^2)$;
\item $\Lc^\eps_\sharp\Phi^\eps=\Lav_\eps\Phi+\mathcal R_\eps$ with $|\mathcal R_\eps|\le C\eps^{2\theta}\mathrm w_p(n)\Pol(u)$; in particular, for $t<\tau^\eps$, $|\mathcal R_\eps(u,m^\eps(t))|\le C\eps^{2\theta-\alpha p}\Pol(u)$;
\item $|\Lav_\eps\Phi(u)-\Lav^N_0\Phi(u)|\le C\eps^{2\theta}\Pol(u)$ and $|\Lav_\eps\Phi(u)|\le C(1+|u|^3+(1+|u|)\|u\|)$.
\end{enumerate}
\end{proposition}

\begin{proof}
(a) Each term of $\Phi^\eps$ is a polynomial in $u$ times a function of $n$ covered by Corollary~\ref{cor:cancel} and Lemma~\ref{lem:R}(iv). For the factors $(u,h)$, the Stokes term is kept in the form $(Au,h)=(u,Ah)$, which is harmless. When stopped at $\tau_R$, the martingale part has an integrable quadratic variation, since the stopped process $\widetilde u(\cdot\wedge\tau_R)$ satisfies $|\widetilde u(t\wedge\tau_R)|\le R$ for all $t\in[0,T]$, and the martingales given by Lemma~\ref{lem:R}(iv) are square integrable on $[0,T]$; hence it is a martingale.

(b) is an easy consequence of the bounds on $\mathsf c,\Psi$ and \eqref{eq:sizes}; for instance $\eps^2|(u,h)(u,\Xi_1)|\le C\eps^{1+2\theta}\mathrm w_{r_2}|u|^2$.

(c) This is Lemma~\ref{lem:remainders}. The last claim follows from Lemma~\ref{lem:stop}.

(d) By \eqref{eq:Y}, $\E_\eta Y-\E_\eta Y_0$ only involves scalar products of $\Psi-\Psi_0$ with $u$, with $\mathsf d+(u,F_h)$, with $A\nN+B(\nN,u_s)$, or with $(\nN\cdot\nabla)u$. Since $|\Psi-\Psi_0|\le C\eps^{2\theta}\mathrm w_{r_1}$ and $\E_\eta[\|n\|_E\mathrm w_{r_1}(n)]<\infty$, the first claim follows. The second follows from \eqref{eq:elem} and the bound on $Y$.
\end{proof}

\subsection{The integrated fast velocity}\label{sec:z}

Recall that $z^\eps(t)=\eps^{-1}\int_0^tm^\eps(s)\dd s$. The process $(u,z^\eps,m^\eps)$ has the formal generator
\[
 \Lc^\eps\Phi+\eps^{-1}D_z\Phi[n],\qquad\Phi=\Phi(u,z,n).
\]
where $\Lc^\eps$ is given by \eqref{eq:generator}.
Since $z^\eps$ is absolutely continuous, Corollary~\ref{cor:cancel} extends without change to polynomials $\Gamma$ in $(u,z)$, with the additional term $\eps^{-1}D_z\Phi[n]$ in $\Lc^\eps_\sharp\Phi$. Fix $h,\phi\in H_{N_0}$.

\begin{proposition}\label{prop:z}
\begin{enumerate}[label=(\alph*)]
\item The process $Z^\eps_\phi(t):=(z^\eps(t),\phi)-\eps(L_1(m^\eps(t)),\phi)+\eps(L_1(m^\eps(0)),\phi)$ is a martingale, and $\E\sup_{t\le T}|z^\eps(t)|^2\le C$.
\item Let $q_\phi(n):=2(n,\phi)(L_1(n),\phi)$, $\chi_\phi:=R_0(q_\phi-\E_\eta q_\phi)$ and
\[
 \Phi^\eps_\phi(z,n):=(z,\phi)^2-2\eps(z,\phi)(L_1(n),\phi)+\eps^2\chi_\phi(n).
\]
Then $\Lc^\eps_\sharp\Phi^\eps_\phi=(Q\phi,\phi)$.
\item Let $\Phi(u,z):=(u,h)(z,\phi)$. There exists $\Phi^\eps=\Phi+\eps\Phi_1+\eps^2\Phi_2$ such that the conclusions (a)--(c) of Proposition~\ref{prop:ptf} hold, with $\mathrm w_p(n)(1+|u|^2)$ replaced by $\mathrm w_p(n)(1+|u|^2+|(z,\phi)|^2)$ in (b), with $\Pol(u)$ replaced by $(1+|(z,\phi)|)\Pol(u)$ in (c), and with $\Lav_\eps\Phi$ replaced by
\begin{gather*}
 \Lav_{\eps,*}\Phi(u,z):=(z,\phi)\Lav_\eps\Phi_h(u)+\E_\eta Y_*(u,\cdot),\\
 Y_*(u,n):=(\sff,h)(L_1(n),\phi)+(n,\phi)(\mathsf c(n)+(u,\Psi(n))),
\end{gather*}
where $\Phi_h(u)=(u,h)$, and $\mathsf c$, $\Psi$ are the coefficients of Section~\ref{sec:ptf} for $\ell=1$. Moreover, if $Y_{*,0}$ denotes $Y_*$ with $\Psi_0$ in place of $\Psi$, then $|\E_\eta Y_*(u,\cdot)-\E_\eta Y_{*,0}(u,\cdot)|\le C\eps^{2\theta}(1+|u|)$.
\end{enumerate}
\end{proposition}

\begin{proof}
Parts (a) and (b) are the classical martingale approximation of the additive functional $z^\eps$ through the Poisson equation $ML_1=\mathrm{Id}$, see~\cite{Bhattacharya1982,KomorowskiLandimOlla2012}. Let us only indicate the main points. For (a), $\Lc^\eps(z,\phi)=\eps^{-1}(n,\phi)$ is compensated by $\eps^{-2}M(-\eps(L_1(n),\phi))=-\eps^{-1}(n,\phi)$. For the bound, the Markov property and \eqref{eq:mix1} give $|\E(m^\eps(r),m^\eps(s))|\le\varrho_0((s-r)/\eps^2)\E_\eta[|n|\mathrm w_{r_1}(n)]$ for $r<s$, hence $\E|z^\eps(t)|^2\le Ct$; since $z^\eps(t)-\eps L_1(m^\eps(t))+\eps L_1(m^\eps(0))$ is an $H$-valued martingale and $\eps^2\E\sup_{t\le T}|L_1(m^\eps(t))|^2\le C\eps^{2-4r_1/q}\le C$ by \eqref{eq:supmoment}, Doob's inequality gives the bound. For (b), the term $-2\eps(z,\phi)(L_1(n),\phi)$ compensates the singular term $2\eps^{-1}(z,\phi)(n,\phi)$ and produces the order-one term $-q_\phi$, which is replaced by its average $-\E_\eta q_\phi=(Q\phi,\phi)$ by the second corrector, see \eqref{eq:k-L1}.

(c) We set
\[
 \Phi_1:=-(u,h)(L_1(n),\phi)+(z,\phi)[\mathsf c(n)+(u,\Psi(n))].
\]
The terms of order $\eps^{-1}$ in $\Lc^\eps_\sharp(\Phi+\eps\Phi_1)$ are $\eps^{-1}(u,h)(n,\phi)-\eps^{-1}(\sff,h)(z,\phi)$, which come from $\Phi$, and $-\eps^{-1}(u,h)(n,\phi)+\eps^{-1}(z,\phi)(\sff,h)$, which come from $\eps^{-2}M(\eps\Phi_1)$ by \eqref{eq:fh} and Corollary~\ref{cor:cancel}; they cancel. The terms of order one which depend on $n$ are $-(z,\phi)Y(u,n)+Y_*(u,n)$, where $Y$ is given by \eqref{eq:Y} with $\ell=1$. By \eqref{eq:fh},
\begin{gather*}
 Y_*-\E_\eta Y_*=z_*(n)+(u,g_*(n)),\\
 z_*:=\mathsf d(L_1,\phi)+(n,\phi)\mathsf c-\E_\eta[\cdots],\qquad g_*:=(L_1,\phi)F_h+(n,\phi)\Psi-\E_\eta[\cdots].
\end{gather*}
These are two-point observables of type (ii) in Lemma~\ref{lem:profiles}, with bilinear maps bounded uniformly in $s$, exactly as $\mathsf d\mathsf c$ and $g_2$ in Section~\ref{sec:ptf}. We set
\[
 \Phi_2:=(z,\phi)[\chi_0+(u,\Xi_1)]-R_0z_*-(u,\PN R_\eps g_*),
\]
with $\chi_0$ and $\Xi_1$ as in Section~\ref{sec:ptf}. The coefficient $\PN R_\eps g_*$ satisfies the bounds \eqref{eq:sizes} of $\Xi_2$. By Corollary~\ref{cor:cancel}, the order-one part of $\Lc^\eps_\sharp(\eps^2\Phi_2)$ is $(z,\phi)(Y-\E_\eta Y)-(Y_*-\E_\eta Y_*)$. The remainder consists of the terms $T_i$ of the proof of Lemma~\ref{lem:remainders} for $\ell=1$, multiplied by $(z,\phi)$, of terms of the same form as $T_1$, $T_3$ and $T_7$, which come from $-(u,h)(L_1(n),\phi)$ and from $\PN R_\eps g_*$, and of the term $\eps(n,\phi)[\chi_0+(u,\Xi_1)]$, which comes from $\eps^{-1}D_z(\eps^2\Phi_2)$ and is bounded by $C\eps^{2\theta}\mathrm w_p(n)(1+|u|)$. The estimates of Lemma~\ref{lem:remainders} and Proposition~\ref{prop:ptf} then give (b) and (c), and (a) is proved in the same way.
\end{proof}

\section{Identification of the effective LU dynamics}\label{sec:identify}

Let us now let $N\to\infty$ in the averaged operators $\Lav^N_0$, and compute the drift and the quadratic variation of the limit equation. The covariance of the noise appears through the Green--Kubo identity \eqref{eq:k-L1}.

\begin{proposition}\label{prop:ident}
There exists $\omega_N\to0$ such that $|\Lav^N_0\Phi(u)-\Lc\Phi(u)|\le\omega_N\Pol(u)$ for $u\in H_N$, where, for $\Phi_h(u)=(u,h)$,
\begin{equation}\label{eq:L-linear}
 \Lc\Phi_h(u)=-(u,Ah)+((u\cdot\nabla)h,u+u_s)-(F(u),h),
\end{equation}
\begin{equation}\label{eq:L-square}
 \Lc\Phi_h^2(u)=2\Phi_h(u)\Lc\Phi_h(u)+\sum_{j\ge0}\mu_j(A\varphi_j+B(\varphi_j,u+u_s),h)^2 .
\end{equation}
Moreover, with $Y_{*,0}$ as in Proposition~\ref{prop:z}, $|\E_\eta Y_{*,0}(u,\cdot)+(Qf_h(u),\phi)|\le\omega_N(1+|u|)$, where $f_h(u)=Ah-P((\nabla h)^T(u+u_s))$.
Moreover, there are functionals $\Lc_J\Phi$, continuous on $H$ and such that $|\Lc_J\Phi(w)|\le C_J(1+|w|^3)$, with $\sup_{N,\eps}\E\int_0^{\tau^\eps}|\Lc\Phi(u)-\Lc_J\Phi(u)|\dd t\le r_J\to0$ as $J\to\infty$.
\end{proposition}

\begin{proof}
\emph{Drift.} Let $N=\infty$, so that $\Psi_0=-B(L_1(n),h)$, and set $v=u+u_s$. By \ref{F4},
\begin{align*}
 -\E_\eta(\sff,\Psi_0)&=\E_\eta(An+B(n,v),B(L_1(n),h))\\
 &=\E_\eta\big(L_1(n),(\nabla h)^T(An+B(n,v))\big)=:-\E_\eta \mathcal B(L_1(n),n),
\end{align*}
where $((\nabla h)^Tw)_i=\sum_k\partial_ih_k\,w_k$, $\mathcal B(w,n)=(w,S_hn)$ and $S_hn=-(\nabla h)^T(An+P((n\cdot\nabla)v))$. The bilinear form $\mathcal B$ is bounded on $E\times E$, so \eqref{eq:k-L1} gives $\E_\eta \mathcal B(L_1(n),n)=-\frac12\mathcal B(k)$, with $\mathcal B$ extended to $\mathcal L_1(E)$.

Let us compute $\mathcal B(k)$ with the eigenvectors of $Q$. We write $k=\sum_is_i\,x_i\otimes y_i$ as in the proof of Lemma~\ref{lem:cov}, so that $\mathcal B(k)=\sum_is_i(x_i,S_hy_i)$. Recall that $Qf=\sum_is_i(x_i,f)y_i$, the series converging in $E$. We may restrict all the sums to the indices $i$ with $s_i>0$. If $Qf=0$, taking the scalar product in $E$ with $y_l$ gives $s_l(x_l,f)=0$, hence $(x_l,f)=0$ whenever $s_l>0$; hence each such $x_i$ is orthogonal to the kernel of $Q$, and $x_i=\sum_j(x_i,\varphi_j)\varphi_j$ in $H$. Since $S_h$ is bounded from $E$ to $H$ and $\sum_is_i<\infty$, we may exchange the sums:
\begin{align*}
 \mathcal B(k)&=\sum_{i,j}s_i(x_i,\varphi_j)(\varphi_j,S_hy_i)=\sum_j\Big(\varphi_j,S_h\sum_is_i(x_i,\varphi_j)y_i\Big)\\
 &=\sum_j(\varphi_j,S_hQ\varphi_j)=\sum_j\mu_j(\varphi_j,S_h\varphi_j).
\end{align*}
The last series converges absolutely by \eqref{eq:noise-smooth}, since $|(\varphi,S_h\varphi)|\le C|\varphi|(|A\varphi|+\|\varphi\|_{L^\infty}\|v\|_{H^1})$. Hence
\begin{align*}
 \E_\eta\mathcal B(L_1(n),n)&=-\frac12\sum_{j\ge0}\mu_j(\varphi_j,S_h\varphi_j)\\
 &=-\frac12\sum_{j\ge0}(B(\sigma_j,A\sigma_j)+B(\sigma_j,B(\sigma_j,v)),h)=(F(u),h),
\end{align*}
where we used $(\varphi,S_h\varphi)=(B(\varphi,A\varphi+B(\varphi,v)),h)$, which follows from \ref{F4}.

Together with $-D\Phi_h[\sfs]=-(u,Ah)-(B(u,v),h)$ and $(B(u,v),h)=-((u\cdot\nabla)h,v)$, this gives \eqref{eq:L-linear}.

\emph{Quadratic variation.} We have $\mathsf c+(u,\Psi_0)=(L_1(n),f_h)$ and $\mathsf d+(u,F_h)=(n,f_h)$, with $f_h=Ah-P((\nabla h)^Tv)$. By \eqref{eq:k-L1},
\begin{equation}\label{eq:qv-green-kubo}
\begin{aligned}
 \E_\eta[(An+B(n,v),h)(AL_1(n)+B(L_1(n),v),h)]&=-\frac12(Qf_h,f_h)\\
 &=-\frac12\sum_{j\ge0}\mu_j(A\varphi_j+B(\varphi_j,v),h)^2 .
\end{aligned}
\end{equation}
Only the symmetric part of the correlation enters here. With \eqref{eq:Y} for $\ell=2$, this gives \eqref{eq:L-square}.

\emph{Cross term.} For $N=\infty$, we have $(\sff,h)=(n,f_h)$ and $\mathsf c+(u,\Psi_0)=(L_1(n),f_h)$, so that $Y_{*,0}=(n,f_h)(L_1(n),\phi)+(n,\phi)(L_1(n),f_h)$. By \eqref{eq:k-L1} and the symmetry of $Q$, $\E_\eta Y_{*,0}=-\frac12(Q\phi,f_h)-\frac12(Qf_h,\phi)=-(Qf_h,\phi)$. The Galerkin error is estimated as below.

\emph{Galerkin error.} Recall that $\Lav^N_0$ differs from $\Lc$ by the replacement of $n$ and $L_1(n)$ by $\PN n$ and $\PN L_1(n)$, and by the outer projection in $\Psi_0$. Since $D(A^{\gamma+1/2})\subset W^{2,\infty}\cap D(A)$, we have, for $w\in E$,
\[
 \|(I-\PN)w\|_{W^{2,\infty}}+|A(I-\PN)w|\le C|A^{\gamma+1/2}(I-\PN)w|\le C\lambda_{N+1}^{-1/2}\|w\|_E .
\] Moreover $|(I-\PN)B(L_1(n),h)|\le\lambda_{N+1}^{-\theta}|A^\theta B(L_1(n),h)|$. Hence, using \eqref{eq:elem}, each difference is bounded by $C\lambda_{N+1}^{-\theta}\mathrm w_p(n)\Pol(u)$ (recall $\theta<1/4$). After averaging in $n$, this gives $\omega_N=C\lambda_{N+1}^{-\theta}$.

\emph{Continuity.} Only the series over $j$ prevent $\Lc\Phi$ from being continuous on $H$: the series in $(F(u),h)$, see \eqref{eq:Fexplicit} and \eqref{eq:L-linear}, and, for $\ell=2$, also the series in \eqref{eq:L-square}. We define $\Lc_J\Phi$ by keeping only the terms $j\le J$ in these series. Let us check that each term is continuous on $H$. Recall that $v=u+u_s$. The map $v\mapsto(B(\varphi_j,v),h)=-((\varphi_j\cdot\nabla)h,v)$ is continuous on $L^2$, with a bound summable by \eqref{eq:noise-smooth}. By \ref{F4}, $(B(\varphi_j,B(\varphi_j,v)),h)=(v,(\varphi_j\cdot\nabla)P((\varphi_j\cdot\nabla)h))$ is continuous on $L^2$ for each $j$. Hence $\Lc_J\Phi$ is continuous on $H$.

Let us now bound the tails $\sum_{j>J}$. For $\ell=1$, the tail of the series in $(F(u),h)$ is bounded by $\sum_{j>J}\mu_j\|\varphi_j\|_{L^\infty}^2\|h\|_{W^{1,\infty}}\|v\|_{H^1}$. For $\ell=2$, the tail is bounded by $2|(u,h)|$ times this quantity, plus the tail of the series in \eqref{eq:L-square}, which is bounded by $2\sum_{j>J}\mu_j(|A\varphi_j|^2|h|^2+\|\varphi_j\|_{L^\infty}^2\|h\|_{W^{1,\infty}}^2|v|^2)$. Since $\E\int_0^{\tau^\eps}(1+|u|)(1+\|u\|)\dd t\le C\,\E\int_0^{\tau^\eps}\Pol(u)\dd t+CT$, the time integrals of these bounds along $u$ on $[0,\tau^\eps]$ tend to $0$ in $L^1(\Omega)$, uniformly in $(N,\eps)$, by Proposition~\ref{prop:energy}; we call $r_J$ this bound.
\end{proof}

\paragraph{Energy balance.} For $u\in V$, \ref{F4} gives $(B(\sigma_j,B(\sigma_j,u)),u)=-|B(\sigma_j,u)|^2$. Hence the part of $-2(F(u),u)$ coming from the transport of $u$ equals $-\sum_j|B(\sigma_j,u)|^2$, and it compensates the corresponding part $\sum_j|B(\sigma_j,u)|^2$ of $\|G(u)\|^2_{L_2(H,H)}$ in It\^o's formula for $|u|^2$. As for any Stratonovich transport noise, the noise acting on $u$ neither creates nor dissipates energy on average; the remaining contributions involve $u_s$ and the additive part $AQ^{1/2}\dd W$. Under \eqref{eq:noise-smooth},
\[
 \|G(u)\|_{L_2(H,H)}^2\le2\sum_j\mu_j|A\varphi_j|^2+2(\sum_j\mu_j\|\varphi_j\|_{L^\infty}^2)|\nabla(u+u_s)|^2\le C(1+\|u\|^2).
\]
Combining \eqref{eq:Fexplicit} and \eqref{eq:Gexplicit}, the limiting equation reads
\begin{align}\label{eq:LUexplicit}
 \dd u+\Big[Au+B(u,u+u_s)&-\frac12\sum_{j\ge0}B(\sigma_j,A\sigma_j)-\frac12\sum_{j\ge0}B(\sigma_j,B(\sigma_j,u+u_s))\Big]\dd t\\
 &=AQ^{1/2}\dd W_t+B(Q^{1/2}\dd W_t,u+u_s).\nonumber
\end{align}
By \eqref{eq:F-decomp}, its transport correction is $\frac12P\nabla\cdot(a\nabla(u+u_s))$, up to the non-local pressure term $\frac12\mathcal R(u+u_s)$.

\begin{remark}[Non-symmetric temporal correlations]
If \eqref{eq:symcorr} is dropped, the proof is unchanged, except for the computation of $\E_\eta\mathcal B(L_1(n),n)$ in Proposition~\ref{prop:ident}. Indeed, \eqref{eq:k-L1} becomes $\E_\eta[L_1(n)\otimes n]=-\frac12(k+k^a)$, where $k$ is still given by \eqref{eq:k} and
\[
 (k^a)^{ij}(x,y)=\int_0^\infty\E[m_i(t,x)m_j(0,y)-m_i(0,x)m_j(t,y)]\dd t
\]
comes from the antisymmetric part of the time-integrated correlation. The quadratic variation \eqref{eq:L-square} is not modified, since only the symmetric part of the correlation enters a quadratic form. The generator \eqref{eq:L-linear} contains the additional term $\frac12\mathcal B(k^a)$, where $\mathcal B$ is the bilinear form of the proof of Proposition~\ref{prop:ident}; it is affine in $u+u_s$. The limit equation is then \eqref{eq:LU} with this additional drift, and it is no longer the It\^o form of the Stratonovich equation \eqref{eq:LU-strat}. Note also that the formal derivation of $u_s$ in Section~\ref{sec:model} uses \eqref{eq:symcorr}. This additional drift is the counterpart of the antisymmetric finite-variation correction to the second level of the limiting rough path obtained in~\cite{DebusscheHofmanova} for an Ornstein--Uhlenbeck fast process; there, it vanishes exactly when the process is reversible.
\end{remark}

\section{Tightness}\label{sec:tight}

We are now ready to prove tightness. We use the Aldous criterion for the finite-dimensional projections of $u$, together with the energy bounds of Proposition~\ref{prop:energy}. Let us first check Aldous' condition for the coordinates of $u$.

\begin{lemma}[Aldous condition for coordinates]\label{lem:aldous}
Let $h\in H_{N_0}$ and $X_t=(u(t),h)$. For every $\eta_0>0$,
\[
 \lim_{\delta_0\to0}\limsup_{\eps\to0}\sup_{N\ge N_0}\sup_{\tau\le T,\ \delta\le\delta_0}\Pp(|X_{(\tau+\delta)\wedge T}-X_\tau|>\eta_0)=0,
\]
where $\tau$ ranges over stopping times.
\end{lemma}

\begin{proof}
Let $\widetilde X_t=(\widetilde u(t),h)$. Since $X=\widetilde X$ on $\{\tau^\eps=T\}$, it suffices to prove the claim for $\widetilde X$, up to the additional term $\Pp(\tau^\eps<T)$, which does not depend on $N$ and tends to $0$. We write $X$ for $\widetilde X$ in the rest of the proof.

Let $\tau_R=\inf\{t:|\widetilde u(t)|\ge R\}\wedge T$, so that $|\widetilde u|\le R$ on $[0,\tau_R]$. Set $\tau'=\tau\wedge\tau_R\wedge\tau^\eps$ and $\tau'_\delta=(\tau+\delta)\wedge\tau_R\wedge\tau^\eps\wedge T$. Then
\[
 \Pp(|X_{(\tau+\delta)\wedge T}-X_\tau|>\eta_0)\le\Pp(\tau_R<T)+\eta_0^{-2}\E|X_{\tau'_\delta}-X_{\tau'}|^2,
\]
and $\Pp(\tau_R<T)\le R^{-2}\E\sup_t|\widetilde u|^2\le C/R^2$ by Proposition~\ref{prop:energy}. Moreover
\[
 \E|X_{\tau'_\delta}-X_{\tau'}|^2=\E[X_{\tau'_\delta}^2-X_{\tau'}^2]-2\E[X_{\tau'}(X_{\tau'_\delta}-X_{\tau'})].
\]
For the first term, Proposition~\ref{prop:ptf} with $\ell=2$ and optional stopping give
\[
 \E[X^2_{\tau'_\delta}-X^2_{\tau'}]\le C\eps\E[\sup_{t\le T}\mathrm w_p(m^\eps(t))](1+R^2)+\E\int_{\tau'}^{\tau'_\delta}|\Lav_\eps\Phi|\dd r+C\eps^{2\theta-\alpha p}\E\int_0^{\tau^\eps}\Pol(u)\dd r .
\]
By \eqref{eq:supmoment} the first term is $O(\eps^{1-2p/q})$, and by (d) and $|\widetilde u|\le R$ on $[\tau',\tau'_\delta]$ the second is bounded by $C_R\delta^{1/2}$.

For the cross term, we use the case $\ell=1$; since $X_{\tau'}$ is $\mathcal F_{\tau'}$-measurable and bounded by $R|h|$, the same bound holds.

We conclude by choosing $R$ large, then $\delta_0$ and $\eps$ small.
\end{proof}

\begin{lemma}[Compactness]\label{lem:compact}
Let $\varsigma>0$. Let $(u^\eps_N)$ be continuous processes in $H$ such that
\[
 \lim_{R\to\infty}\limsup_{\eps\to0}\sup_N[\Pp(\sup_t|u^\eps_N(t)|>R)+\Pp(\int_0^T\|u^\eps_N\|^2\dd t>R)]=0,
\]
and such that, for every $k$, the real processes $(u^\eps_N,e_k)$ satisfy the Aldous condition~\cite[Section~16]{Billingsley1999}, in the asymptotic form of Lemma~\ref{lem:aldous}. Let
\[
 \mathcal X_T=C([0,T];D(A^{-\varsigma}))\cap L^2(0,T;H).
\]
Then for every $\kappa>0$ there exist a compact set $\mathcal K_\kappa\subset\mathcal X_T$ and $\eps_\kappa>0$ such that $\Pp(u^\eps_N\in\mathcal K_\kappa)\ge1-\kappa$ for $\eps<\eps_\kappa$ and all $N$ large enough.
\end{lemma}

\begin{proof}
The proof is classical; we only note that all the bounds are needed only for $\eps$ small. Let $\Pi_K$ denote the orthogonal projection onto $H_K=\mathrm{span}(e_1,\dots,e_K)$; here $K$ is a truncation parameter, independent of the Galerkin level $N$.

For each $K$, Aldous' criterion~\cite[Section~16]{Billingsley1999} applied to $\Pi_Ku^\eps_N$ gives the required compact sets in the Skorokhod space $D([0,T];H_K)$; for continuous processes, the $J_1$ modulus controls the uniform one, so these sets can be taken compact in $C([0,T];H_K)$. Moreover $\sup_t|(I-\Pi_K)u^\eps_N(t)|_{D(A^{-\varsigma})}\le\lambda_{K+1}^{-\varsigma}\sup_t|u^\eps_N(t)|$ is small in probability, uniformly in $N$ and for $\eps$ small. A standard approximation argument then gives the statement in $C([0,T];D(A^{-\varsigma}))$.

Finally, a set bounded in $L^2(0,T;V)$ and relatively compact in $C([0,T];D(A^{-\varsigma}))$ is relatively compact in $L^2(0,T;H)$, by the interpolation inequality $|w|\le C|A^{-\varsigma}w|^{1/(1+2\varsigma)}\|w\|^{2\varsigma/(1+2\varsigma)}$ and H\"older's inequality in time.
\end{proof}

Let us check that the assumptions of Lemma~\ref{lem:compact} hold for $u=u^\eps_N$. Indeed, $u$ and $\widetilde u$ coincide on $[0,\tau^\eps]$, so that
\[
 \Pp(\sup_t|u(t)|>R)+\Pp(\int_0^T\|u\|^2\dd t>R)\le2\Pp(\tau^\eps<T)+\frac{\E\sup_{t\le\tau^\eps}|u|^2}{R^2}+\frac{\E\int_0^{\tau^\eps}\|u\|^2\dd t}R,
\]
which is small for $R$ large and $\eps$ small by Proposition~\ref{prop:energy}; the Aldous condition is Lemma~\ref{lem:aldous}. Note that the stopped process $\widetilde u$ itself is not used here, since $\|\widetilde u(t)\|=\|u(\tau^\eps)\|$ for $t>\tau^\eps$ is not controlled.

Hence, by Lemma~\ref{lem:compact} and Prokhorov's theorem, for every $N(\eps)\to\infty$, the family $\{\Law(u^\eps_{N(\eps)})\}$ is tight in $\mathcal X_T$ as $\eps\to0$.

\begin{lemma}[Tightness of the integrated fast velocity]\label{lem:ztight}
The family $\{\Law(z^\eps)\}$ is tight in $\mathcal Z_T$ as $\eps\to0$.
\end{lemma}

\begin{proof}
Let $\phi\in H_{N_0}$. The proof of Lemma~\ref{lem:aldous} applies to $(z^\eps,\phi)$, with Proposition~\ref{prop:z}(a)--(b) in place of Proposition~\ref{prop:ptf}, and with the stopping time $\inf\{t:|z^\eps(t)|\ge R\}\wedge T$ in place of $\tau_R$. It is even simpler, since $\Lc^\eps_\sharp\Phi^\eps_\phi=(Q\phi,\phi)$ has no remainder, and since $\eps\E\sup_{t\le T}\mathrm w_{p'}(m^\eps(t))\to0$ for $p'\le2p$ by \eqref{eq:supmoment}. Hence the real processes $(z^\eps,e_k)$ satisfy the Aldous condition. Since $\E\sup_{t\le T}|z^\eps(t)|^2\le C$ by Proposition~\ref{prop:z}(a), the first part of the proof of Lemma~\ref{lem:compact} gives the result.
\end{proof}

Hence the family $\{\Law(u^\eps_{N(\eps)},z^\eps)\}$ is tight in $\mathcal X_T\times\mathcal Z_T$, which proves the first assertion of Theorem~\ref{thm:main}(ii).
\section{Passage to the limit and proof of Theorem~\ref{thm:main}}\label{sec:passage}

We now conclude the proof of Theorem~\ref{thm:main}. Let us first explain the strategy. We pass to the limit in the martingale identities given by the perturbed test functions: first for the observables of $u$, which gives a martingale problem for the limit $u$, and then for the observables involving $z$, which identifies the noise. Thanks to the stopping times, all the error terms are controlled by Proposition~\ref{prop:energy}, uniformly in $N$ and $\eps$.

Let $\eps_k\downarrow0$, $N_k=N(\eps_k)$, $u^k=u^{\eps_k}_{N_k}$, $z^k=z^{\eps_k}$, $\widetilde u^k=u^k(\cdot\wedge\tau^{\eps_k})$ and $\widetilde z^k=z^k(\cdot\wedge\tau^{\eps_k})$. Since $\Pp(\widetilde u^k\ne u^k)+\Pp(\widetilde z^k\ne z^k)\le2\Pp(\tau^{\eps_k}<T)\to0$, the laws of $(\widetilde u^k,\widetilde z^k)$ have the same limit points as those of $(u^k,z^k)$, and they are tight along the sequence by Section~\ref{sec:tight}. By Skorokhod's theorem, along a subsequence, there exist, on another probability space, processes with the same laws as $(\widetilde u^k,\widetilde z^k)$ which converge almost surely in $\mathcal X_T\times\mathcal Z_T$ to a limit $(u,Z)$; we keep the notation $(\widetilde u^k,\widetilde z^k)$ for them. Only these processes are transported: all the computations below which involve $m^k$ or $\tau^{\eps_k}$ are performed on the original probability space, and they lead to bounds and identities which only depend on the law of $(\widetilde u^k,\widetilde z^k)$.

\emph{Regularity of the limit.} Let us first check the regularity of $u$. The map $w\mapsto\sup_t|w(t)|^2$ is lower semicontinuous on $\mathcal X_T$, so Fatou's lemma and Proposition~\ref{prop:energy} give $\E\sup_t|u(t)|^2\le\liminf_k\E\sup_{t\le\tau^{\eps_k}}|u^k(t)|^2\le C$.

For the $V$-norm, we cannot use $\widetilde u^k$ directly, since $\|\widetilde u^k(t)\|=\|u^k(\tau^{\eps_k})\|$ for $t>\tau^{\eps_k}$ is not controlled. Let $\|w\|^2:=\sup_K\sum_{k\le K}\lambda_k(w,e_k)^2\in[0,\infty]$. For $M>0$, the map $w\mapsto\int_0^T(\|w(t)\|^2\wedge M)\dd t$ is lower semicontinuous on $L^2(0,T;H)$, and
\[
 \E\int_0^T(\|\widetilde u^k\|^2\wedge M)\dd t\le\E\int_0^{\tau^{\eps_k}}\|u^k\|^2\dd t+MT\,\Pp(\tau^{\eps_k}<T)\le C+MT\,\Pp(\tau^{\eps_k}<T).
\]
Fatou's lemma gives $\E\int_0^T(\|u\|^2\wedge M)\dd t\le C$, and letting $M\to\infty$ yields $u\in L^2(\Omega\times(0,T);V)$.

\emph{Martingale problem.} Fix $N_0$, $h\in H_{N_0}$, $\ell\in\{1,2\}$, $0\le s<t\le T$ and a bounded continuous functional $\Xi$ of the path of $(u,z)$ on $[0,s]$. Let $\tau^k_R:=\inf\{r:|\widetilde u^k(r)|\ge R\}\wedge T$ and, for $r\in[0,T]$, $r_R:=r\wedge\tau^k_R\wedge\tau^{\eps_k}$. By Proposition~\ref{prop:ptf}(a) and optional stopping,
\[
 \E\Big[\Xi(\widetilde u^k,\widetilde z^k)\Big(\Phi^{\eps_k}(\widetilde u^k_{t_R},m^k_{t_R})-\Phi^{\eps_k}(\widetilde u^k_{s_R},m^k_{s_R})-\int_{s_R}^{t_R}\Lc^{\eps_k}_\sharp\Phi^{\eps_k}\dd r\Big)\Big]=0 .
\]
By Proposition~\ref{prop:ptf}(b), \eqref{eq:supmoment} and the quartic bound of Proposition~\ref{prop:energy} at the stopping times $t_R$ and $s_R$, replacing $\Phi^{\eps_k}$ by $\Phi$ costs at most $C\eps_k(\E\sup_r\mathrm w_{2p}(m^k(r)))^{1/2}(1+\E|\widetilde u^k_{t_R}|^4)^{1/2}\le C\eps_k^{1-2p/q}$, and similarly at time $s$.

By Propositions~\ref{prop:ptf}(c)--(d) and~\ref{prop:ident}, we may replace $\Lc^{\eps_k}_\sharp\Phi^{\eps_k}$ by $\Lc\Phi(\widetilde u^k)$, up to an error $C(\eps_k^{2\theta-\alpha p}+\eps_k^{2\theta}+\omega_{N_k})$, and then $\Lc\Phi$ by $\Lc_J\Phi$, up to $r_J$. All these errors are uniform in $R$.

We now remove the stopping times. Since $|\Phi(w)|+|\Lc_J\Phi(w)|\le C_J(1+|w|^3)$, and since $\E|\widetilde u^k(\sigma)|^4\le C$ for every stopping time $\sigma\le T$ and $\E\int_0^T|\widetilde u^k(r)|^4\dd r\le CT$ by Proposition~\ref{prop:energy}, the random variables $\Phi(\widetilde u^k_{t_R})$, $\Phi(\widetilde u^k_{s_R})$ and $\int_0^T|\Lc_J\Phi(\widetilde u^k_r)|\dd r$ are uniformly integrable, uniformly in $k$ and $R$. Moreover, the times $r_R$ and $r$ differ only on $\{\tau^k_R<T\}\cup\{\tau^{\eps_k}<T\}$, and $\Pp(\tau^k_R<T)\le C/R^2$ by Proposition~\ref{prop:energy}, while $\Pp(\tau^{\eps_k}<T)\to0$. Recall that $\widetilde u^k$ is constant after $\tau^{\eps_k}$. Letting $R\to\infty$, we have thus shown that
\[
 \Big|\E\Big[\Xi(\widetilde u^k,\widetilde z^k)\Big(\Phi(\widetilde u^k_t)-\Phi(\widetilde u^k_s)-\int_s^t\Lc_J\Phi(\widetilde u^k_r)\dd r\Big)\Big]\Big|\le Cr_J+o(1)\qquad\text{as }k\to\infty .
\]
This quantity only depends on the law of $(\widetilde u^k,\widetilde z^k)$. Passing to the limit $k\to\infty$ and then $J\to\infty$,
\[
 \E[\Xi(u,Z)(\Phi(u(t))-\Phi(u(s))-\int_s^t\Lc\Phi(u(r))\dd r)]=0 .
\]
Thus $M_h(t)=(u(t),h)-(u_0,h)-\int_0^t\Lc\Phi_h(u(r))\dd r$ is a continuous martingale for the filtration generated by $(u,Z)$, and the case $\ell=2$ together with \eqref{eq:L-square} gives, by polarization,
\[
 \langle M_h,M_g\rangle_t=\int_0^t(G(u(r))^*h,G(u(r))^*g)\dd r =\int_0^t(Qf_h(u(r)),f_g(u(r)))\dd r .
\]

\emph{The integrated fast velocity.} Let $\phi\in H_{N_0}$, and let $\Phi$ be one of the functions $(z,\phi)$, $(z,\phi)^2$ and $(u,h)(z,\phi)$, with the perturbed test functions of Proposition~\ref{prop:z}. The corresponding limit operators are $\Lc\Phi=0$, $\Lc\Phi=(Q\phi,\phi)$ and $\Lc\Phi(u,z)=(z,\phi)\Lc\Phi_h(u)-(Qf_h(u),\phi)$, by Proposition~\ref{prop:z} and Proposition~\ref{prop:ident}.

We argue as above, with one more stopping time $\vartheta^k_R:=\inf\{r:|\widetilde z^k(r)|_{D(A^{-\varsigma})}\ge R\}\wedge T$, so that $|(\widetilde z^k,\phi)|\le R|A^\varsigma\phi|$ on $[0,\vartheta^k_R]$. The stopping time $\tau^k_R$ is removed as above, by uniform integrability, but we keep $\vartheta^k_R$, since we do not have a fourth moment bound on $z^\eps$.

Since $Z$ is continuous with values in $D(A^{-\varsigma})$, for almost every $R$ the hitting time $\vartheta_R$ is almost surely a continuity point of the corresponding functional at $Z$, so that $\vartheta^k_R\to\vartheta_R:=\inf\{r:|Z(r)|_{D(A^{-\varsigma})}\ge R\}\wedge T$ almost surely. Passing to the limit $k\to\infty$, we obtain that $\Phi(u,Z)(\cdot\wedge\vartheta_R)-\int_0^{\cdot\wedge\vartheta_R}\Lc\Phi(u,Z)\dd r$ is a martingale. Letting $R\to\infty$, the processes
\[
 Z_\phi:=(Z,\phi),\qquad Z_\phi^2-t(Q\phi,\phi),\qquad(u,h)Z_\phi-\int_0^\cdot\big(Z_\phi\Lc\Phi_h(u)-(Qf_h(u),\phi)\big)\dd r
\]
are continuous local martingales for the filtration generated by $(u,Z)$. By It\^o's formula for the product $(u,h)Z_\phi$, this means that
\[
 \langle Z_\phi\rangle_t=t(Q\phi,\phi),\qquad\langle M_h,Z_\phi\rangle_t=-\int_0^t(Qf_h(u(r)),\phi)\dd r .
\]
By L\'evy's characterization~\cite{DaPratoZabczyk2014}, $Z$ is a $Q$-Wiener process for this filtration.

\emph{Identification of the noise.} Let $N_h:=M_h+\int_0^\cdot(f_h(u(r)),\dd Z_r)$. The stochastic integral is well defined, since $|Q^{1/2}f_h(u)|\le C(1+|u|)$. By the above brackets and a density argument,
\[
 \langle N_h\rangle_t=\int_0^t\big((Qf_h,f_h)-2(Qf_h,f_h)+(Qf_h,f_h)\big)(u(r))\dd r=0,
\]
so that $N_h=0$. On an enlarged probability space, let $W$ be a cylindrical Wiener process on $H$ such that $Q^{1/2}W=-Z$. Since $(w,f_h(u))=(Aw+B(w,u+u_s),h)$ for $w\in E$ by \ref{F4}, we obtain
\[
 M_h(t)=-\int_0^t(A\dd Z_r+B(\dd Z_r,u+u_s),h)=\int_0^t(G(u(r))\dd W_r,h)
\]
for all $h$ in $\bigcup_{N_0}H_{N_0}$, hence for all admissible $h$. By \eqref{eq:L-linear}, this is the weak form of \eqref{eq:LUexplicit}. Therefore every limit point of $\Law(u^\eps_{N(\eps)},z^\eps)$ is the law of a pair $(u,Z)$ as in Theorem~\ref{thm:main}(ii), which proves (ii). Part (iii) follows from (ii) by the beginning of Section~\ref{sec:energy}, and (i) was proved in Section~\ref{sec:existence}.
\begin{proposition}[Pathwise uniqueness for the limit equation in dimension two]\label{prop:LU2D}
Let $d=2$. Let $u_1,u_2$ be two martingale solutions of \eqref{eq:LU} defined on the same stochastic basis, with the same cylindrical Wiener process $W$ and the same initial datum. Then $u_1=u_2$ almost surely.
\end{proposition}

\begin{proof}
Let $\beta_j=(W,\varphi_j)$, which are independent Brownian motions, and $w=u_1-u_2$.

Since $G(u)$ vanishes on the kernel of $Q$, $G(u)\dd W=\sum_jG(u)\varphi_j\dd\beta_j$. Since $G(u_1)-G(u_2)=B(Q^{1/2}\cdot,w)$ and $F(u_1)-F(u_2)=-\frac12\sum_jB(\sigma_j,B(\sigma_j,w))$, the process $w$ satisfies, in $V'$,
\[
 \dd w+[Aw+B(w,u_1+u_s)+B(u_2,w)-\frac12\sum_jB(\sigma_j,B(\sigma_j,w))]\dd t=\sum_jB(\sigma_j,w)\dd\beta_j .
\]
In dimension two, the drift belongs to $L^2(0,T;V')$ almost surely, by Ladyzhenskaya's inequality and \eqref{eq:noise-smooth}, and the noise coefficient belongs to $L^2(0,T;L_2(H,H))$, since $\sum_j|B(\sigma_j,w)|^2\le\sum_j\mu_j\|\varphi_j\|_{L^\infty}^2\|w\|^2$. Therefore It\^o's formula in the Gelfand triple~\cite{KrylovRozovskii1981} applies and gives
\begin{align*}
 \dd|w|^2&=[-2\|w\|^2-2(B(w,u_1+u_s),w)-2(B(u_2,w),w)\\
 &\qquad+\sum_j(B(\sigma_j,B(\sigma_j,w)),w)+\sum_j|B(\sigma_j,w)|^2]\dd t+2\sum_j(B(\sigma_j,w),w)\dd\beta_j .
\end{align*}
By \ref{F4}, $(B(u_2,w),w)=0$, $(B(\sigma_j,w),w)=0$ and $(B(\sigma_j,B(\sigma_j,w)),w)=-|B(\sigma_j,w)|^2$. Hence the It\^o correction cancels exactly, the martingale term vanishes, and
\begin{align*}
 \frac{\dd}{\dd t}|w|^2+2\|w\|^2&=-2(B(w,u_1+u_s),w)=2(B(w,w),u_1)-2(B(w,u_s),w)\\
 &\le\|w\|^2+C(\|u_1\|^2+\|\nabla u_s\|_{L^\infty})|w|^2,
\end{align*}
pathwise, by Ladyzhenskaya's inequality. Since $\|u_1\|^2\in L^1(0,T)$ almost surely, the result follows from Gronwall's lemma.
\end{proof}

\emph{Proof of Theorem~\ref{thm:main}(iv).} Thanks to Proposition~\ref{prop:LU2D} and the Yamada--Watanabe theorem in infinite dimensions~\cite{RocknerSchmulandZhang2008}, \eqref{eq:LU} has a unique strong solution: $u$ is a measurable functional of $W$, and the joint law of $(u,W)$ is unique. Since $Z=-Q^{1/2}W$, the law of the pair $(u,Z)$ in (ii) and (iii) is unique. Hence all the limit points in (ii) and (iii) coincide, and, together with the tightness, $(u^\eps,z^\eps)$ and $(u^\eps_{N(\eps)},z^\eps)$ converge in law to $(u,Z)$ as $\eps\to0$.

\section{Spatially homogeneous fast advection}\label{sec:homogeneous}

In a bounded domain, the limit noise cannot be spatially homogeneous. Indeed, $m$ vanishes on $\partial S$, so that the variance tensor $a(x)=k(x,x)$ vanishes on $\partial S$. If $a$ were constant, it would vanish identically. Since $\operatorname{tr}a(x)=\sum_j\mu_j|\varphi_j(x)|^2$ with $\mu_j>0$, all the $\varphi_j$ would vanish, that is $Q=0$, and the limit equation would be deterministic. In order to study a spatially homogeneous noise, we therefore consider the torus. On the torus $\T^d=\R^d/(2\pi\mathbb Z)^d$, with $H$ and $V$ replaced by spaces of mean-zero divergence-free fields, the whole analysis applies with simplifications: there is no boundary layer, and one may take $\eps=0$ in $R_\eps$. We now assume that the statistics of $m$ are invariant under translations, so that $k(x,y)=k_0(x-y)$ with
\begin{gather*}
 k_0(z)=\sum_{\zeta\in\mathbb Z^d_0}\hat k_0(\zeta)e^{i\zeta\cdot z},\\
 \hat k_0(\zeta)\ \text{Hermitian, non-negative},\qquad\hat k_0(-\zeta)=\overline{\hat k_0(\zeta)},\qquad\hat k_0(\zeta)\zeta=0,
\end{gather*}
where $\mathbb Z^d_0=\mathbb Z^d\setminus\{0\}$ and $\sum_\zeta|\zeta|^4\operatorname{tr}\hat k_0(\zeta)<\infty$, which is the analogue of \eqref{eq:noise-smooth} in this setting. Then $a=k_0(0)=\sum_\zeta\hat k_0(\zeta)$ is constant and real symmetric, and $u_s=0$. Note also that $\xi\cdot\hat k_0(\zeta)\xi$ is real for $\xi\in\R^d$. The computations below only use these properties. For $\xi\in\mathbb Z^d_0$, we write $\Pi_\xi=I-\xi\otimes\xi/|\xi|^2$ for the Leray projector on the Fourier mode $\xi$, and $\Pi_0=0$.

\emph{The drift associated with the additive noise vanishes.} Since $\sum_j\sigma_j(x)\otimes\sigma_j(y)=k_0(x-y)$ and each $\sigma_j$ is divergence free, $\sum_j(\sigma_j\cdot\nabla)A\sigma_j=-\nu\nabla\cdot(\sum_j\sigma_j\otimes\Delta\sigma_j)$. But $\sum_j\sigma_j(x)\otimes\Delta\sigma_j(x)=\Delta k_0(0)$ does not depend on $x$. Hence $\sum_jB(\sigma_j,A\sigma_j)=0$.

\emph{The transport correction is an explicit Fourier multiplier.} Let $\mathcal T(v):=\sum_jB(\sigma_j,B(\sigma_j,v))$. Since $\mathcal T$ is quadratic in $(\sigma_j)$, it only depends on $\sum_j\sigma_j\otimes\sigma_j=k$, and we may compute it with the pairs of complex modes $\alpha e^{i\zeta\cdot x}$ and $\bar\alpha e^{-i\zeta\cdot x}$, $\alpha\otimes\bar\alpha$ summing to $\hat k_0(\zeta)$. Let $v=\hat ve^{i\xi\cdot x}$ with $\hat v\perp\xi$. The inner transport gives $P((\bar\alpha e^{-i\zeta\cdot x}\cdot\nabla)v)=i(\bar\alpha\cdot\xi)\Pi_{\xi-\zeta}\hat ve^{i(\xi-\zeta)\cdot x}$. The outer transport and projection then give $-(\alpha\cdot(\xi-\zeta))(\bar\alpha\cdot\xi)\Pi_\xi\Pi_{\xi-\zeta}\hat ve^{i\xi\cdot x}$. Summing over $\alpha$ and $\zeta$, and using $\hat k_0(\zeta)\zeta=0$, hence $\hat k_0(\zeta)^T\zeta=\overline{\hat k_0(\zeta)\zeta}=0$, so that $(\xi-\zeta)\cdot\hat k_0(\zeta)\xi=\xi\cdot\hat k_0(\zeta)\xi$, we obtain
\begin{equation}\label{eq:T-symbol}
 \mathcal T(v)=-[\sum_{\zeta\in\mathbb Z^d_0}(\xi\cdot\hat k_0(\zeta)\xi)\,\Pi_\xi\Pi_{\xi-\zeta}]\hat v\,e^{i\xi\cdot x}.
\end{equation}
Since $P\nabla\cdot(a\nabla v)=-(\xi\cdot a\xi)\hat ve^{i\xi\cdot x}=-\sum_\zeta(\xi\cdot\hat k_0(\zeta)\xi)\hat ve^{i\xi\cdot x}$, the pressure term \eqref{eq:F-decomp} is the Fourier multiplier
\begin{equation}\label{eq:R-symbol}
 \mathcal R(v)=-[\sum_{\zeta\ne\xi}(\xi\cdot\hat k_0(\zeta)\xi)\,\Pi_\xi\frac{(\xi-\zeta)\otimes(\xi-\zeta)}{|\xi-\zeta|^2}\Pi_\xi]\hat v\,e^{i\xi\cdot x}.
\end{equation}
Its symbol is a non-positive symmetric matrix on $\xi^\perp$. It vanishes only if $\hat k_0(\zeta)\xi=0$ for every $\zeta$ which is not parallel to $\xi$; in general $\mathcal R\ne0$, and it is not a multiple of the Laplacian. The limiting equation reads
\[
 \dd u+[Au+B(u,u)-\frac12\mathcal T(u)]\dd t=AQ^{1/2}\dd W_t+B(Q^{1/2}\dd W_t,u),\qquad-\frac12\mathcal T=-\frac12a:\nabla^2+\frac12\mathcal R .
\]
Moreover $-\frac12\mathcal T$ is a non-negative operator, since $(\mathcal T(v),v)=-\sum_j|B(\sigma_j,v)|^2$ by \ref{F4}; this is the energy balance of Section~\ref{sec:identify}.

\emph{Isotropic noise concentrated on high modes.} Let $\rho\in C^1_c((1,2))$ be non-negative and non-zero, and for $L\ge1$ let
\begin{gather*}
 \hat k_0^L(\zeta)=L^{-d}\rho(|\zeta|/L)\,\Pi_\zeta,\qquad a^L=a_0^LI,\\
 a_0^L=\frac{d-1}d\,L^{-d}\sum_\zeta\rho(|\zeta|/L)\longrightarrow a_0:=\frac{d-1}d\int_{\R^d}\rho(|y|)\dd y .
\end{gather*}
Then $\xi\cdot\hat k_0^L(\zeta)\xi=L^{-d}\rho(|\zeta|/L)|\xi|^2\sin^2\angle(\zeta,\xi)$, and $|\Pi_{\xi-\zeta}-\Pi_\zeta|\le C|\xi|/L$ on the support. Hence, for fixed $\xi$,
\[
 \hat{\mathcal T}^L(\xi)=-|\xi|^2L^{-d}\sum_\zeta\rho(|\zeta|/L)\sin^2\angle(\zeta,\xi)\,\Pi_\xi\Pi_\zeta\Pi_\xi+O(|\xi|^3/L).
\]
The sum is a Riemann sum, with mesh $1/L$, of the $C^1$ function $y\mapsto\rho(|y|)\sin^2\angle(y,\xi)\Pi_\xi\Pi_y\Pi_\xi$, supported in $\{1\le|y|\le2\}$. It therefore equals $\int_{\R^d}\rho(|y|)\sin^2\angle(y,\xi)\Pi_\xi\Pi_y\Pi_\xi\dd y+O(1/L)$. In polar coordinates, this integral is $(\int\rho(|y|)\dd y)\E[\sin^2\angle(\hat y,\xi)\Pi_\xi\Pi_{\hat y}\Pi_\xi]$, with $\hat y$ uniform on the sphere. By symmetry, the expectation is a multiple $\beta_d\Pi_\xi$; taking the trace and using $\E[\hat y_d^2]=\frac1d$ and $\E[\hat y_d^4]=\frac3{d(d+2)}$ gives $\beta_d=\frac{d^2-3}{d(d+2)}$. Altogether,
\[
 \hat{\mathcal T}^L(\xi)=-c_d\,a_0|\xi|^2\Pi_\xi+O((|\xi|^2+|\xi|^3)/L),\qquad\hat{\mathcal R}^L(\xi)=-(1-c_d)\,a_0|\xi|^2\Pi_\xi+O((|\xi|^2+|\xi|^3)/L),
\]
with $c_d=\frac{d}{d-1}\beta_d=\frac{d^2-3}{(d+2)(d-1)}$, that is $c_2=\frac14$ and $c_3=\frac35$. Hence $\mathcal T^L\to c_da_0\Delta$ and $\mathcal R^L\to(1-c_d)a_0\Delta$ mode by mode, with rate $O(1/L)$. In dimension three, the factor $c_3=\frac35$ also appears in the high-mode limit of the transport noise considered by Flandoli and Luo~\cite{FlandoliLuo2021}. With the normalization of~\cite{DebusscheHofmanova2026Fluct}, for which $a_0=\frac{d-1}dq(0)$ in this limit, the resulting enhanced viscosity $\frac12c_da_0$ equals $\frac{q(0)}{16}$ for $d=2$ and $\frac{q(0)}5$ for $d=3$, in agreement with their Lemma~4.3. At fixed covariance, however, \eqref{eq:T-symbol}--\eqref{eq:R-symbol} are the exact operators, and the pressure correction is genuinely non-local.

\end{document}